\documentclass[10pt,a4paper,reqno,twoside]{amsart}

\usepackage{amsmath,amssymb,amsthm}
\usepackage[english]{babel}
\usepackage[latin1]{inputenc}
\usepackage{enumerate}
\usepackage{mathtools}
\usepackage{color}
\usepackage{fancyhdr}
\usepackage{enumitem}

\def\@begintheorem#1#2{\bgroup{\sc #1 \ #2.}\it \ignorespace}
\def\@opargbegintheorem#1#2#3{ \bgroup{\sc #1\ #2 \ (#3).\\}\it \ignorespace}
\def\@endtheorem{\egroup}

\newtheorem{thm}{Theorem}[section]
 
 \newtheorem{lem}[thm]{Lemma}
 \newtheorem{prop}[thm]{Proposition}
 \theoremstyle{definition}
 \newtheorem{defn}[thm]{Definition}
 \theoremstyle{remark}
 \newtheorem{rem}[thm]{Remark}
\theoremstyle{definition}
\newtheorem{hypothesis}{Hypothesis}
 \newtheorem{ex}{Example}
 \numberwithin{equation}{section}

\newcommand{\eps}{\varepsilon}

\newcommand{\R}{\mathbb R}

\renewcommand{\phi}{\varphi}

\newcommand{\N}{\mathbb N}
\newcommand{\PI}{\mathcal P}
\newcommand{\supp}{\mbox{supp\,}}
\newcommand{\lin}{{\mathrm{\,lin\,}}}
\newcommand{\loc}{{\mathrm{\,loc\,}}}

\newcommand{\etab}{\bar{\eta}}
\newcommand{\epsb}{\bar{\eps}}
\newcommand{\argchi}{\big(tR^{-1}, xR^{-\frac{1}{\etab}}\big)}

\newenvironment{proofof}[1]{\smallskip\noindent{\textbf{Proof~of~#1.}}%
	\hspace{1pt}}{\hspace{-5pt}{\nobreak\quad\nobreak\hfill\nobreak%
		$\square$\vspace{2pt}\par}\smallskip\goodbreak}

\title[Critical Non-linearity for some evolution equations with Fujita-type critical exponent]{Critical non-linearity for some evolution equations with Fujita-type critical exponent }
\date{\today}
\author[G. Girardi]{Giovanni Girardi}

\address{Giovanni Girardi, Department of Industrial Engineering and Mathematical Sciences, Polytechnic University of Marche, Via Brecce Bianche 12 - 60131 ANCONA - ITALY}
\email{g.girardi@univpm.it
}

\begin{document}
\maketitle
\begin{abstract}
We consider the Cauchy problem for a class of non-linear evolution equations in the form
\begin{equation*}
\label{eq}
L(\partial_t,\partial_x) u=F(\partial_t^\ell u), \quad (t,x)\in [0,\infty)\times \mathbb{R}^n;
\end{equation*}
here, $L(\partial_t,\partial_x)$ is a linear partial differential operator with constant coefficients, of order $m\geq 1$ with respect to the time variable $t$, and $\ell$ is a natural number satisfying $0\leq \ell\leq m-1$.  For several different choices of $L$, many authors have investigated the existence of global (in time) solutions to this problem when $F(s)=|s|^p$ is a power non-linearity, looking for a \textit{critical exponent} $p_c>1$ such that  global small data solutions exist in the supercritical case $p>p_c$, whereas no global weak solutions exist, under suitable sign assumptions on the data, in the subcritical case $1<p<p_c$.\\
In the present paper we consider a more general non-linear term in the form $F(s)=|s|^p\mu(|s|)$; for a large class of models, we provide an integral condition on $\mu$ which allows to distinguish more precisely the region of existence of a global (in time) small data solution from that in which the problem admits no global (in time) weak solutions, refining the existing results about the critical exponents for power type non-linearities.
\end{abstract}



\section{Introduction}
We consider differential operators in the form
\begin{equation}
	\label{eq:Ldef}
	L(\partial_t,\partial_x)=\partial_t^m+ \sum_{j=0}^{m-1} P_{j}(\partial_x)\partial_t^j,
\end{equation}
where $m$ is a positive integer number, and for all $j=0,1,\dots, m-1$ it holds
$$P_{j}(\partial_x)=\sum_{|\alpha|= 0}^{d_j} c_{j,\alpha} \partial_x^\alpha, $$ 
for some $d_j\in \N$ and $c_{j,\alpha}\in \R$; for ease of presentation we also fix $d_m=0$ and $P_m(\partial_x)=1$.

In recent years, many authors have investigated the existence of global (in time) solutions to the Cauchy problem associated to the non-linear equation 
\begin{equation}
\label{eq:CPnon-linear_classic}
L(\partial_t,\partial_x) u = F(\partial_t^\ell u),  \qquad (t,x)\in [0,\infty)\times \R^n, \quad \ell=0,\dots, m-1
\end{equation}
with power non-linearity $F(s)=|s|^p$ for $p> 1$, under suitable assumptions on the initial data.\\
The presence of a non-linear term in the equation make in general not applicable the boot-strap argument which allows to prolong local in time solutions to well-posed linear problems to large time, proving the existence of global (in time) solutions. This issue motivates the study of conditions on the differential operator $L$, on the non-linear function $F$ and on the initial data, under which global (in time) solutions to the non-linear Cauchy problem associated to \eqref{eq:CPnon-linear_classic} exist, or do not exist. 
\subsection{The critical exponent}
In the pioneering work \cite{Fujita}, the author H. Fujita considered the Cauchy problem for the heat equation ($L=\partial_t-\Delta$) with power non-linearity $|u|^p$ and he provided a threshold power $p_c=1+2/n$ such that global (in time) solutions exist for $p>p_c$, if the initial data are sufficiently small, and do not exist if $p\in (1, p_c] $ for any positive initial data, even under smallness assumptions. This threshold power $p_c$ is nowadays known as \textit{critical exponent}. More in general, we say that $p_c>1$ is a critical exponent for problem \eqref{eq:CPnon-linear_classic} with power non-linearity $F(s)=|s|^p$ if a global (in time) solution exists for $p>p_c$  for sufficiently small initial data in a suitable space, whereas no global weak solutions exist if $1 < p < p_c$ under a suitable sign condition on the data, even if initial data are arbitrarily small; the exponent $p=p_c$ may belong to the existence or non-existence range of global solutions, according to the problem considered.\\
The existing literature has shown that the presence of terms of different nature in the equation, and their interactions, can deeply influence the critical exponent; an overview about existing results about wave-type models and $\sigma$-evolution equations, with different kind of dissipations is provided, for instance, in \cite{DAG2024} and the reference therein. \\
For some models, the critical exponent $p_c$ is related to the scaling properties of the differential operator $L$; in this case, the critical power is often referred as \textit{Fujita-type critical exponent}, as in this manuscript,  and the application of a test function method allows in general to prove the non-existence of global weak solutions for $p\in (1,p_c]$ (we refer the reader to \cite{DAL2013,DAL2003, MitPohoz98,MitPohoz99,MitPohoz01} for details about the test function method, and to  \cite{Lucente} for a review about wave-type problems with Fujita-type critical exponents); whereas, for $p>p_c$ the existence of global (in time) small data solutions to the Cauchy problem associated to \eqref{eq:CPnon-linear_classic} can sometimes be proved
treating the non-linearity as a perturbation (for small initial data), employing appropriate long time $L^p-L^q$ decay estimates for the solution to the corresponding linear problem (see the approach used in Section \ref{sec:existence}), proving the sharpness of the non-existence results. 

In \cite{DAF2021} the authors provide a constructive method which allows to directly compute the critical exponent of non-existence for Cauchy problem \eqref{eq:CPnon-linear_classic} with power non-linearity $F(s)=|s|^p$, when the operator $L$ is in the form 
\[L= \partial_t^m+\sum_{j=0}^{m-1}a_j(-\Delta)^{r_j}\partial_t^j,\]
for some $a_j\in\R$, simply knowing the powers $r_j$ which are integer or fractional numbers.
With minor modifications the approach used in \cite{DAF2021} can be applied to obtain a constructive method for the critical exponent $p_c$ of non-existence for non-linear problem \eqref{eq:CPnon-linear_classic} with $F(s)=|s|^p$, and $L$ in the more general form \eqref{eq:Ldef}; for later use, we provide it hereafter. Let us consider the set 
\begin{equation}
\label{eq:J} 
J=\{ j\in \{0,\dots, m\}: P_j(\partial_x)\neq 0 \},
\end{equation}
and, for all $j\in J$, we define
\[ r_j:=\min\{ |\alpha|: c_{j,\alpha}\neq 0\};\]
then, we introduce the functions $g:\eta\in [0,\infty]\to \R$ such that 
\begin{equation*}
g(\eta):=\min_{\substack{j\in J}}\{(j-\ell)\eta+r_j\},
\end{equation*}
and $h:[0,\infty]\to \R$ such that
\begin{equation*}
h(\eta):=1+\frac{g(\eta)}{(n+\eta-g(\eta))_+},
\end{equation*}
where $a_+=\max(a,0)$ for $a\in \R$, and we set $1/0=\infty$. Then, we define 
\begin{equation} 
\label{eq:pcritical}
p_c:=\max_{\eta\in [0,\infty]}h(\eta).
\end{equation}
If $p_c\in (1,\infty)$, applying a classical test function method it is possible to prove that no nontrivial global weak solutions to \eqref{eq:CPnon-linear_classic} exist for  $1<p\leq p_c$, under suitable sign assumptions on the initial data. The authors in \cite{DAF2021} provide a detailed explanation about how the representation \eqref{eq:pcritical} for the critical exponent $p_c$ can be derived and how to compute it; in particular, they show that if $p_c\in (1,\infty)$ then there exists $\bar{\eta}\in (0,\infty)$ such that $p_c=h(\bar{\eta})$. \\
Whereas, if $p_c=1$ the application of a test function method does not allow in general to get any result about existence or non-existence of global (in time) solutions to \eqref{eq:CPnon-linear_classic}. In this case, it is possible for some models to prove that a global (in time) small data solution to \eqref{eq:CPnon-linear_classic} exists for all $p>1$, provided that the solution to the corresponding linear problem decays fast enough for large time (see later Proposition \ref{thm:existence2}), as it happens for instance for the damped $\sigma$-evolution equation with mass (see Seciton \ref{sec:KG}).\\
In order to make clearer the given procedure, here we provide an example of application. 
\begin{ex}
\label{example:pc}
Consider the Cauchy problem for the following evolution equation with structural damping 
\begin{equation}
\label{eq:structurally_damped}
\partial_t^2 u +(-\Delta)^{\delta} \partial_tu+(-\Delta)^{\sigma}u =F(u), \quad F(s)=|s|^p
\end{equation}
with $0\leq \delta<\sigma$ natural numbers; for this model we find $m=2$ and $J=\{0,1,2\}$ with $r_0=2\sigma$, $r_1=2\delta$, $r_2=0$; moreover, it holds
\begin{align*}
g(\eta)&=\min\{2\eta, \eta+2\delta, 2\sigma\};
\end{align*}
in order to apply formula \eqref{eq:pcritical} we need to distinguish three cases. \\
If $\delta=0$, we get 
\[g(\eta)=\begin{cases}
\eta \quad &\text{ if } \eta\in [0,2\sigma], \\
2\sigma \quad &\text{ if } \eta\in (2\sigma,\infty];
\end{cases}\]
then, $\bar{\eta}=2\sigma$ and
\[ p_c=h(2\sigma)=1+\frac{2\sigma}{n}.\]
If $0<\delta<\sigma/2$, we find
\[g(\eta)=\begin{cases}
\eta \quad &\text{ if } \eta\in [0,2\delta], \\
\eta+2\delta  \quad &\text{ if } \eta\in (2\delta,2(\sigma-\delta)],\\
2\sigma \quad &\text{ if } \eta\in (2(\sigma-\delta),\infty];
\end{cases}\]
thus, in this case we get $\bar{\eta}=2(\sigma-\delta)$ and
\[ p_c=h(2(\sigma-\delta))=1+\frac{2\sigma}{n-2\delta},\]
for all $n>2\delta$.
Finally, if $\delta\geq \sigma/2>0$ we obtain 
\[g(\eta)= \begin{cases}
2\eta \quad &\text{ if } \eta\in [0,\sigma], \\
2\sigma \quad &\text{ if } \eta\in (\sigma,\infty];
\end{cases}\]
then, $\bar{\eta}=\sigma$ and
\[ p_c=h(\sigma)=1+\frac{2\sigma}{n-\sigma}.\]
provided that $n>\sigma$.\\
No global (in time) weak solutions to \eqref{eq:structurally_damped} exists if $p$ satisfies $1<p\leq p_c$;
on the other hand, the existence of a global (in time) small data solution can be proved, in low space dimension, for all $p>p_c$ provided that $2\delta\leq \sigma$ (see \cite{DAE2017}) or $2\delta>\sigma>1$ (see \cite{DAE22}); additional details about this model are provided in Section \ref{sec:structurally}.
\end{ex}
\noindent Other examples of evolution equations with Fujita-type critical exponents, expressible through formula \eqref{eq:pcritical}, are provided and discussed in Section \ref{sec:examples}. \\
Fixing $\sigma=1$ and  $\delta=0$ in \eqref{eq:structurally_damped} we find the non-linear classical damped wave equation 
\begin{equation}
\label{eq:damped_wave_intro}
\partial_t^2 u + \partial_tu-\Delta u =F(u), \quad F(s)=|s|^p
\end{equation}
whose critical exponent is the Fujita exponent $p_c=1+2/n$, the same found in \cite{Fujita} for the semilinear heat equation (see \cite{ IkeTanizawa2005, Matsumura, TodoYord2001, Zhang}); indeed, the classical damping term $\partial_t u$ produces a \textit{diffusion phenomenon}, i.e., the solution to the linear problem associated to \eqref{eq:damped_wave_intro} behaves asymptotically like the solution to an heat equation with suitable initial data (see \cite{MarcatiNish2003,Nishihara2003,YangMilani}). In particular, no global (in time) weak solutions exist for all $p\leq 1+2/n$.

\subsection{A new critical non-linearity}
In \cite{EGR} we have shown that the known results about the existence of global (in time) small data solutions to the Cauchy problem associated to \eqref{eq:damped_wave_intro} can be refined in space dimension $n=1, 2$, considering a more general non-linear term in the form $F(s)=|s|^{1+2/n}\mu(|s|)$ where $\mu$ is a modulus of continuity: under suitable assumptions on $\mu$ the Cauchy problem associated to \eqref{eq:damped_wave_intro} admits a global (in time) small data solution if $\mu$ satisfies the integral condition 
\begin{equation}\label{eq:muint}
	\int_0^{c_0} \frac{\mu(\tau)}\tau\,d\tau <\infty,
\end{equation}
for some $c_0>0$ arbitrarily small; whereas, no global weak solutions exist if \eqref{eq:muint} is not fulfilled for any $c_0>0$, under a suitable sign condition on the initial data.\\
Similarly, in \cite{DAG2024} we considered the Cauchy problem associated to \eqref{eq:structurally_damped}, with $2\delta>\sigma>1$; we showed that global (in time) small data solutions still exist if one repleace the power non-linearity $|u|^p$ by the non-linear term $F(u)=|u|^{1+2\sigma/(n-\sigma)}\mu(|u|)$ where $\mu$ is, for instance, a modulus of continuity satisfying the integral condition \eqref{eq:muint}, for some $c_0>0$ arbitrarily small.

Motivated by the results obtained in \cite{DAG2024} and \cite{EGR}, in the present paper we consider the Cauchy problem associated to equation \eqref{eq:CPnon-linear_classic}, with $F(s)=|s|^p\mu(|s|)$, and we show that for a large class of differential operators $L$ the integral condition \eqref{eq:muint} on $\mu$ provides a critical assumption on the non-linearity $F$ which allows to distinguish more precisely the region of existence of a global (in time) small data solution from that in which no global solutions exist, under suitable sign assumptions on the initial data, refining the existing results about the critical exponents for power type non-linearities. \\
More in detail, suppose that $p_c\in (1,\infty)$ is the critical exponent for equation \eqref{eq:CPnon-linear_classic} with power non-linearity $F(s)=|s|^p$, and consider the Cauchy problem associated to the non-linear model
\begin{equation}
\label{eq:CPnon-linear_mu}
L(\partial_t,\partial_x) u = |\partial_t^\ell u|^{p_c}\mu(|\partial_t^\ell u|),  \qquad (t,x)\in [0,\infty)\times \R^n,
\end{equation}
where $\mu:[0,\infty)\to [0,\infty)$ satisfies the following assumptions:
\begin{enumerate}[label=(\roman*), itemsep=0pt, topsep=0pt]
\item $\mu$ is non-decreasing;
\item $ \mu$ is bounded on $[0,\epsb]$, for some $\epsb>0$;
\item the non-linear function $F: \R \to \R$, $F(s):= |s|^{p_c}\mu(|s|)$ is convex;
\item the function $F: \R \to \R$, $F(s)= |s|^{p_c}\mu(|s|)$ verifies $F(0)=0$ and
\begin{equation}
\label{eq:Fcontraction_intro}
|F(y)-F(z)|\leq C\,|y-z|\,(|y|^{p_c-1}+|z|^{p_c-1})\,\mu(|y|+|z|),
\end{equation}
for all $|y|\leq\epsb$ and $|z|\leq\epsb$, for some constant $C>0$ independent of $y$ and $z$.
\end{enumerate}
For a large class of models we show that, problem \eqref{eq:CPnon-linear_mu} admits a unique global small data solution if $\mu$ satisfies the integral condition \eqref{eq:muint} for some $c_0>0$ arbitrarily small; whereas, no global weak solutions exist if \eqref{eq:muint} is not fulfilled, under a suitable sign condition on the initial data.\\
On the other hand, if $p_c=1$ then condition \eqref{eq:muint} is in general not necessary for the existence of nontrivial global (in time) solutions to \eqref{eq:CPnon-linear_mu}. \\
The paper is organized as follows: in Section \ref{sec:non-existence} we apply an improved test function method introduced in \cite{IkeSob}, proving that if $p_c\in (1,\infty)$ and $\mu$ does not satisfy \eqref{eq:muint} then, no nontrivial global weak solutions to \eqref{eq:CPnon-linear_mu} exist (see Definition \ref{def:weaksolution}) under a suitable sign assumption on the initial data; in Section \ref{sec:existence} we show that the integral condition \eqref{eq:muint} on $\mu$ is also sufficient to guarantee the global (in time) existence of small data solutions to \eqref{eq:CPnon-linear_mu}, provided that the solution to the corresponding linear problem satisfies suitable long time decay estimates. For the sake of completness, in Section \ref{sec:existence-p=1} we also consider the case $p_c=1$. In Section \ref{sec:examples} we give some examples in which the obtained results can be applied. Finally, in Section \ref{sec:conclusions} we provide some additional remarks about the obtained results and we propose some open problems.

\subsubsection{Some remarks on the function $\mu=\mu(s)$}
In the following, we provide some remarks useful to better understand which functions $\mu=\mu(s)$ one could consider.
\begin{rem}
\label{rem:mu}
Notice that if $\mu$ is continuous then condition \eqref{eq:muint} implies $\mu(0)=0$.
\end{rem}
\begin{rem}
If $\mu$ satisfies
\begin{equation*}
\label{eq:mu_assumptions_existence}
\mu \in \mathcal C^1, \quad  0\leq\tau\mu'(\tau)\leq\mu(\tau) \text{ for } \tau\in[0,\epsb],
\end{equation*} 
then assumption~\eqref{eq:Fcontraction_intro} holds true for $|y|\leq \epsb$ and $|z|\leq\epsb$, provided that $p_c\geq 1$. Indeed,
\[ F(y)-F(z)=\int_0^1 \partial_\rho F(z+\rho(y-z))\,d\rho = (y-z)\int_0^1 F'(z+\rho(y-z))\,d\rho, \]
so that
\[\begin{split}
|F(y)-F(z)|
    & \leq |y-z| \int_0^1 |F'(z+\rho(y-z))|\,d\rho \\
    & \leq p_c|y-z|\,\int_0^1 |z+\rho(y-z)|^{p_c-1}\,\mu(|z+\rho(y-z)|)\,d\rho\\
    & \leq p_c|y-z|\,(|y|+|z|)^{p_c-1}\,\mu(|y|+|z|).
\end{split} \]
\end{rem}

\begin{ex}
\label{example:mu}
We recall that a \textit{modulus of continuity} is a function $\mu\in\mathcal C([0,\infty))$, increasing and concave, that satisfies $\mu(0)=0$. A modulus of continuity, clearly meets conditions (i) and (ii). \\ As an example of function $\mu$ one can consider a modulus of continuity defined for a sufficiently small~$\tau$ by $(-\log\tau)^{-\gamma}$ or, more in general, by
\[ \mu(s)=\Big(-\log\tau\Big)^{-1}\Big(\log(-\log\tau)\Big)^{-1}\cdots \Big(\log^{[k]}(-\log \tau)\Big)^{-\gamma},\quad k \in \N,\]
where $\log^{[k]}$ denotes the composition of $k$ $\log$, and for large $\tau$ in such a way that assumptions (i)-(iv) are satisfied. With this choice of $\mu$, condition \eqref{eq:muint} is satisfied if, and only if, $\gamma>1$.
\end{ex}
%
\section{Non-existence of global solutions} 
\label{sec:non-existence}
Let $L$ be a differential operator as in \eqref{eq:Ldef} and $\ell$  a natural number between $0$ and $m-1$; suppose that the critical power $p_c$ defined by \eqref{eq:pcritical} belongs to $(1,\infty)$. In this section we prove a non-existence result for the non-linear Cauchy problem
\begin{equation}
\label{eq:CP_non-linear_non-existence}
\begin{cases}
L(\partial_t,\partial_x)u=F(\partial_t^\ell u), &\quad t>0,\,x\in \R^n, \\
\partial_t^j u(0,x)=u_j(x), &\quad j=0,\dots, m-1,
\end{cases}
\end{equation}
where $F:\R\to \R$ is defined by $F(s)=|s|^{p_c}\mu(|s|)$, with $\mu:[0,\infty)\to [0,\infty)$ satisfying assumptions (i)-(iii).
%

%
%
\begin{defn}
\label{def:weaksolution}
Let $L$ be as in \eqref{eq:Ldef} and $\ell\in \{0,\dots, m-1\}$. Assume that the initial data in \eqref{eq:CP_non-linear_non-existence} satisfy
\[ u_j=0 \;\text{ if }\, j=0,\dots,\ell-1 \;\text{ and }\, u_j\in L^1_{\text{loc}}(\R^n)\; \text{ if }\, j\geq \ell.\]
Let us fix $T\in (0,\infty]$; we say that $u\in W^{\ell,p_c}_{\text{loc}}\big([0,T), L^{p_c}_{\text{loc}}(\R^n)\big)$ is a weak solution to \eqref{eq:CP_non-linear_non-existence} if it satisfies the initial conditions $\partial_t^j u(0,\cdot)=u_j$ for every $j\in \{0,\dots, m-1\}$, and for any function $\psi\in C_c^\infty([0,T)\times \R^n)$ satisfying $\psi\equiv 1$ in a neighborhood of $0$, it holds
\begin{align*}
\notag
	\int_0^T \int_{\R^n}F(&\partial_t^\ell u(t,x))\psi(t,x)\,dx\,dt \\
	%
	&=\sum_{j=0}^m(-1)^{j-\ell}\int_0^T \int_{\R^n}\partial_t^\ell u(t,x) \tilde{P}_{j}(\partial_x)\partial_t^{(j-\ell)}\psi(t,x)\,dx\,dt\\
&-\sum_{j=\ell}^{m-1}\int_{\R^n}u_j(x)\tilde{P}_{j+1}(\partial_x)\psi(0,x)\,dx;
\end{align*}
here, for any $j<0$, $\partial_t^{(j)}\psi$ is defined as the primitive of $\partial_t^{(j+1)}\psi$, i.e.
\[ \partial_t^{(j)}\psi(t,x)=-\int_t^T\partial_t^{(j+1)}\psi(\tau,x)\,d\tau,\]
which is still compactly supported; moreover, 
\[ \tilde{P}_{j}(\partial_x)=\sum_{|\alpha|=0}^{d_j} (-1)^{|\alpha|}c_{j,\alpha}\partial_x^\alpha.\]
\end{defn}

\begin{thm}
\label{thm:non-existenceCauchy}
Let $I=\{j\geq \ell: c_{j+1,0}\neq 0\}$; we assume that $u_j=0$ for any $j\leq \ell-1$ and $u_j\in L^1(\R^n)$ for any $j=\ell,\dots,m-1$, together with the sign condition
\begin{equation}
\label{eq:sign_assumption}
 \sum_{j\in I} c_{j+1,0}\int_{\R^n}u_j(x)\,dx>0.
\end{equation}
Suppose that the exponent $p_c$ defined by formula \eqref{eq:pcritical} belongs to $(1,\infty)$ and $\mu: [0,\infty)\to [0,\infty)$ verifies assumptions (i)-(iii).\\
If there exists a global-in-time weak solution $u\in W^{\ell,p_c}_{\text{loc}}\big([0,\infty),L^{p_c}_{\text{loc}}(\R^n)\big)$ to \eqref{eq:CP_non-linear_non-existence}, according to Definition \ref{def:weaksolution}, then,
\begin{equation}
\label{eq:blowupcondition}
\int_{0}^{c_0} \frac{\mu(s)}{s}\,ds<\infty,
\end{equation}
for some $c_0>0$ sufficiently small.
\end{thm}
\begin{rem}
We note that $I\neq \emptyset$, indeed $j=m-1$ belongs to $I$.
\end{rem}
In the proof of Theorem \ref{thm:non-existenceCauchy} we will employ the following generalized Jensen inequality.
\begin{lem}
	\label{Lemma:Jensen}
	Let $F$ be a convex function on $\R$. Let $\alpha=\alpha(x)$ be defined and non-negative almost everywhere on $\Omega$, such that $\alpha$ is positive in a set of positive measure. Then, it holds
	\begin{equation*}
	\label{eq:Jensen}
	 F\bigg(\frac{\int_\Omega v(x)\alpha(x)\,dx}{\int_\Omega \alpha(x)\,dx}\bigg)\leq \frac{\int_\Omega F(v(x))\alpha(x)\,dx}{\int_\Omega \alpha(x)\,dx} 
	\end{equation*}
	for all non-negative functions $v$ provided that all the integral terms are meaningful.
\end{lem}

\begin{proofof}{Theorem \ref{thm:existence}}
Let $\etab\in [0,+\infty)$ such that $p_c=h(\bar{\eta})$. We consider %
\[ \chi(s):=\begin{cases}
1 & \quad \text{ if } s\in [0,1/2],\\
0 & \quad \text{ if } s\in [1,+\infty),\\
\text{decreasing } & \quad \text{ if } s\in (1/2,1),
\end{cases}
%
\]
where $\chi\in C^\infty ([0,+\infty))$. We define the cut-off function on $[0,+\infty)\times \R^n$
\[\psi(t,x):=(\chi(t+|x|^{\etab}))^q,\]
where $q:=\max_{j\in J}\{d_j+(j-\ell)_+\}p'_c$, with $p'_c$ satisfying $1/p_c+1/p'_c=1$. \\
Furthermore, for all $R>0$ we set
\[\psi_R(t,x)=\psi\argchi;
\]
we note that 
\[ \supp(\psi_R)\subset Q_R:= [0,R]\times B_{R^\frac{1}{\etab}}(0). \]
%
%
Let us suppose by contradiction that $\mu=\mu(s)$ does not satisfy the integral condition \eqref{eq:blowupcondition} for any $c_0>0$ and problem \eqref{eq:CP_non-linear_non-existence} admits a global (in time) weak solution in the sense of Definition \ref{def:weaksolution}; we define
\[ I_R=\int_0^{+\infty}\int_{\R^n} 
F(\partial_t^\ell u(t,x))\psi_R(t,x)\,dx\,dt.\] Then, it holds
\begin{align*}
I_R =\sum_{j=0}^m(-1)^{j-\ell}\int_0^{+\infty} \int_{\R^n}\partial_t^\ell u(t,x) & \tilde P_j(\partial_x)\partial_t^{(j-\ell)}\psi_R(t,x)\,dx\,dt \\ &-\sum_{j=\ell}^{m-1}\int_{\R^n}u_j(x)\tilde P_{j+1}(\partial_x)\psi_R(0,x)\,dx.
\end{align*}
It is easy to check that for all $j\in J$ it holds
\begin{align*}
\tilde P_{j}(\partial_x)\partial_t^{(j-\ell)}&\psi_R(t,x)\\
&= R^{-(j-\ell)}\sum_{|\alpha|=r_j}^{d_j} R^{-\frac{|\alpha|}{\etab}} (-1)^{|\alpha|}c_{j,\alpha} (\partial_x^\alpha\partial_t^{(j-\ell)} \psi)\argchi;
\end{align*}
in particular, we may employ the boundeness of $\chi$ and all its derivatives to estimate 
\begin{equation*}
|\partial_x^\alpha \partial_t^{(j-\ell)}\psi|\lesssim \psi ^{\frac{q-j+\ell-|\alpha|}{q}},
\end{equation*}
if $j\geq \ell$,  being $q\geq j-\ell+|\alpha|$ for all $r_j\leq |\alpha|\leq d_j$.\\
On the other hand, since $\supp \psi\subset Q_1$ and $\psi$ is non-increasing, we may estimate
\begin{equation*}
|\partial_x^\alpha \partial_t^{(j-\ell)}\psi|\lesssim \psi^{\frac{q-|\alpha|}{q}},
\end{equation*}
if $j<\ell$. 
Finally, being $\psi\leq 1$ we may conclude
\begin{align*}
I_R \lesssim \sum_{j\in J} R^{-(j-\ell)-\frac{r_j}{\etab}} & \int_0^{+\infty} \int_{\R^n}|\partial_t^\ell u(t,x)|  \psi_R(t,x) ^{\frac{q-d_j-(j-\ell)_+}{q}}\,dx\,dt \\ &-\sum_{j=\ell}^{m-1}\int_{\R^n}u_j(x)\tilde{P}_{j+1}(\partial_x)\psi_R(0,x)\,dx.
\end{align*}
for all $R>1$.
Let us introduce the auxiliary functions
\[ \lambda(t,x)=(t+|x|^{\etab})^\beta(1+t+|x|^{\etab})^{-1},\quad  \lambda_R(t,x)=\lambda\argchi,\]
for some $\beta<\min\{1,p_c-1\}$; since $F$ is a convex function, for all $j\in J$ we can apply the generalized Jensen inequality given in Lemma \ref{Lemma:Jensen} with 
\[v(t,x):= |\partial_t^\ell u(t,x)|\lambda_R(t,x)^{\frac{1}{p_c-1}}\psi_R(t,x)^{\frac{d_j+(j-\ell)_+}{q(p_c-1)}}\]
and 
\[ \alpha(t,x):=\lambda_R(t,x)^{-\frac{1}{p_c-1}}\psi_R(t,x)^{1-\frac{p_c(d_j+(j-\ell)_+)}{q(p_c-1)}}\]
to get 
\begin{align*}
&F\left(\frac{\int_0^{+\infty} \int_{\R^n}|\partial_t^\ell u(t,x)|\psi_R(t,x)^{\frac{q-d_j-(j-\ell)_+}{q}}\,dx\,dt}{\int_0^{+\infty} \int_{\R^n}\alpha(t,x)\,dx\,dt}\right)\\
&\hspace{40pt}=F\left(\frac{\int_0^{+\infty} \int_{\R^n}v(t,x)\alpha(t,x)\,dx\,dt}{\int_0^{+\infty} \int_{\R^n}\alpha(t,x)\,dx\,dt}\right) \\
&\hspace{40pt}\lesssim \frac{\int_0^{+\infty}\int_{\R^n}F(v(t,x))\alpha(t,x)\,dx\,dt}{\int_0^{+\infty} \int_{\R^n}\alpha(t,x)\,dx\,dt}.
\end{align*}
Since $\mu$ is increasing, $\lambda_R\leq 1$ and $\psi_R\leq 1$ we may estimate
\[ F(v(t,x))\alpha(t,x)\leq F(\partial_t^\ell u(t,x))\lambda_R(t,x)\psi_R(t,x)\leq F(\partial_t^\ell u(t,x))\psi_R(t,x);\]
on the other hand, for all $j\in J$ we may estimate
\begin{align*}
&\int_0^{+\infty} \int_{\R^n}\alpha(t,x)\,dx\,dt\\ & \hspace{10pt}= \int_0^{+\infty} \int_{\R^n}\lambda(tR^{-1}, x R^{-\frac{1}{\etab}})^{-\frac{1}{p_c-1}}\psi(tR^{-1}, x R^{-\frac{1}{\etab}})^{1-\frac{p_c(d_j+(j-\ell)_+)}{q(p_c-1)}}\,dx\,dt\\
& \hspace{10pt}= R^{1+\frac{n}{\eta}} \int_{Q_1} \lambda(\tau,y)^{-\frac{1}{p_c-1}}\psi(\tau, y)^{1-\frac{p_c(d_j+(j-\ell)_+)}{q(p_c-1)}}\,dy\,d\tau \\
& \hspace{10pt}= C_{j,d_j} R^{1+\frac{n}{\etab}},
\end{align*}
where \[C_{j,d_j} = \int_{Q_1} (\tau+|y|^{\etab})^{-\frac{\beta}{p_c-1}}(1+\tau+|y|^{\etab})^{\frac{1}{p_c-1}}\psi(\tau, y)^{1-\frac{p_c(d_j+(j-\ell)_+)}{q(p_c-1)}}\,dy\,d\tau;\]
note that, for all $j\in J$ the quantity $C_{j,d_j}$ is a finite positive number; indeed, our choice of $q$ allows to estimate 
\[ 0< (\tau+|y|^{\etab})^{-\frac{\beta}{p_c-1}}(1+\tau+|y|^{\etab})^{\frac{1}{p_c-1}}\psi(\tau, y)^{1-\frac{p_c(d_j+(j-\ell)_+)}{q(p_c-1)}}\lesssim \tau^{-\frac{\beta}{p_c-1}},\]
for all $(\tau, y) \in Q_1$; in particular, our choice of $\beta$ guarantees that the right-hand side has finite integral over $Q_1$.

All of these considerations allow to conclude that there exists $C>0$ such that the following estimate holds true:
\begin{align*}
I_R &+ \sum_{j=\ell}^{m-1}\int_{\R^n}u_j(x)\tilde P_{j+1}(\partial_x)\psi_R(0,x)\,dx \\ & \lesssim  \sum_{j\in J} R^{-(j-\ell)-\frac{r_j}{\etab}+1+\frac{n}{\etab}} F^{-1}\bigg(C R^{-1-\frac{n}{\etab}}\int_0^{+\infty}\int_{\R^n} F(\partial_t^\ell u)\lambda_R\psi_R\,dx\,dt\bigg).
\end{align*}
For all $r>0$ we define 
\[ y(r):=\int_0^{+\infty}\int_{\R^n}F(\partial_t^\ell u(t,x))\lambda_r(t,x)\psi_r(t,x)\,dx\,dt,\]
and 
\[ Y(R):= \int_0^R y(r)r^{-1}\,dr. \]
On the one hand, being $y(R)=Y'(R)R$ we find
\begin{equation} 
\label{eq:preliminary_1}
\begin{aligned}
I_R  + \sum_{j=\ell}^{m-1}&\int_{\R^n}u_j(x)\tilde P_{j+1}(\partial_x)\psi_R(0,x)\,dx\\ & \lesssim \sum_{j\in J} R^{-(j-\ell)-\frac{r_j}{\etab}+1+\frac{n}{\etab}} F^{-1}\big(C R^{-\frac{n}{\etab}}Y'(R)\big), 
\end{aligned}
\end{equation}
for all $R>1$. On the other hand, using that $\psi_r\leq \psi_R$ for all $r\in [0,R]$ and the change of variable $s=(t+|x|^{\bar{\eta}})$ we get
\begin{equation}
\label{eq:preliminary_2}
\begin{aligned}
Y(R)&= \int_0^{+\infty}\int_{\R^n}F(\partial_t^\ell u(t,x))\bigg( \int_0^R \lambda_r(t,x)\psi_r(t,x) r^{-1}\,dr \bigg)\,dx\,dt \\
& \lesssim  \int_0^{+\infty}\int_{\R^n}F(\partial_t^\ell u(t,x))\psi_R(t,x) \bigg( \int_{\frac{t+|x|^{\etab}}{R}}^{+\infty} \frac{s^{\beta-1}}{1+s}\,ds \bigg)\,dx\,dt \\
& \lesssim  \int_0^{+\infty}\int_{\R^n}F(\partial_t^\ell u(t,x))\psi_R(t,x)\,dx\,dt = I_R.
\end{aligned}
\end{equation}
Collecting \eqref{eq:preliminary_1} and \eqref{eq:preliminary_2} we obtain
\begin{align*}
Y(R)\lesssim  \sum_{j\in J} &R^{-(j-\ell)-\frac{r_j}{\etab}+1+\frac{n}{\etab}} F^{-1}\big(C R^{\frac{n}{\etab}}Y'(R)\big)\\&-\sum_{j=\ell}^{m-1}\int_{\R^n}u_j(x)\tilde P_{j+1}(\partial_x)\psi_R(0,x)\,dx,
\end{align*}
for all $R>1$.
We note that for all $j\geq \ell$ such that $j\notin I$ the support of $\tilde P_{j+1}(\partial_x)\psi_R(0,x)$ is a subset of $B_R(0)\setminus B_{R/2}(0)$; as consequence, since $u_j\in L^1(\R^n)$ it holds
\[ \lim_{R\to +\infty}\int_{\R^n}u_j(x)\tilde P_{j+1}(\partial_x)\psi_R(0,x)\,dx=0;\]
then, employing also the definition of $g(\etab)$ and the sign condition \eqref{eq:sign_assumption}, we may estimate for $R>R_0\geq 1$ sufficiently large
\begin{align*} Y(R)&\leq C_0 R^{-\frac{g(\etab)}{\etab}+1+\frac{n}{\bar{\etab}}} F^{-1}\big(C R^{-\frac{n}{\etab}}Y'(R)\big),\\
&= C_0 R^{\frac{n+\etab}{p_c \etab}}F^{-1}\big(C R^{-\frac{n}{\etab}}Y'(R)\big)
\end{align*}
for a suitable constant $C_0>0$, and then
\begin{align*}
 \frac{1}{C R}\mu\bigg(\frac{Y(R_0)}{C_0 R^\frac{n+\etab}{p_c\etab}}\bigg)\leq \frac{Y'(R)C_0^{p_c}}{Y(R)^{p_c}}.
\end{align*}
Integrating from $R_0$ to $R$ we obtain the existence of a positive constant $c_1$ such that
\begin{align*} 
\int_{R_0}^R \frac{1}{s}\mu\bigg(\frac{Y(R_0)}{C_0 s^\frac{n+\etab}{p_c\etab}}\bigg)\,ds &= c_1 \int_{R^{-\frac{n+\etab}{p_c\etab}}}^{R_0^{-\frac{n+\etab}{p_c\etab}}}\frac{\mu(s)}{s}\,ds \\ &\lesssim \left[Y(R)^{1-p_c}\right]_{R_0}^R\lesssim Y(R_0)^{1-p_c}.
\end{align*}
A contradiction follows for any choice of $R_0>1$ taking the limit for $R\rightarrow \infty$; indeed, the right-hand side is a finite quantity, whereas the integral term on the left-hand side cannot be bounded since $\mu$ does not satisfy assumption \eqref{eq:blowupcondition} for any $c_0>0$. 
The proof of the statement follows.
\end{proofof}

\section{Existence of global solutions}\label{sec:existence}
\subsection{The case $p_c>1$}
Let $L$ be a partial differential operator as in  \eqref{eq:Ldef} and  $\ell$ a natural number in  $\{0,\dots, m-1\}$. In this section we provide some global (in time) existence results for the following non-linear Cauchy  problem
\begin{equation}
	\label{eq:CPnon-linear}
	\begin{cases}
		L(\partial_t,\partial_x) u = F(\partial_t^\ell u), \\
		\partial_t^j u(0,x)=0, \quad  \text{ for all } j=0,\dots, m-2,\\
		\partial_t^{m-1} u(0,x)= f(x),
	\end{cases}
\end{equation}
where $F$ is a non-linear term defined by $F(s)=|s|^{p_c}\mu(|s|)$, with $p_c> 1$ given by formula \eqref{eq:pcritical} and $\mu:[0,\infty)\to [0,\infty)$ a function satisfying properties (i), (ii) and (iv).\\
Given $K=K(t,x)$ the fundamental solution for the partial differential operator $L$, for every $f$ we can write the solution to the corresponding linear Cauchy problem 
\begin{equation}
	\label{eq:CPlinear}
	\begin{cases}
		L(\partial_t,\partial_x) u = 0, \\
		\partial_t^j u(0,x)=0, \quad  \text{ for all } j=0,\dots, m-2,\\
		\partial_t^{m-1} u(0,x)= f(x),
	\end{cases}
\end{equation}
as 
\[ u^\lin(t,x)=K(t,\cdot)\ast_{(x)} f;\]
moreover, due to Duhamel's principle the solution to \eqref{eq:CPnon-linear} can be written as 
\begin{equation} 
\label{eq:duhamel}
u(t,x)=u^\lin (t,\cdot)+\int_0^t K(t-s,\cdot)\ast_{(x)}F(\partial_t^\ell u(s,\cdot))\,ds. \quad 
\end{equation}
The representation \eqref{eq:duhamel} allows to prove the following result, employing a standard contraction argument together with suitable $L^p-L^q$ estimates for the linear Cauchy problem \eqref{eq:CPlinear}.
\begin{thm}
\label{thm:existence}
Suppose that $p_c>1$ and there exists $r\geq 1$ such that for every $q\in [p_c,\infty]$ and $f\in L^1\cap L^r$ the solution to the linear problem \eqref{eq:CPlinear} satisfies
\begin{equation}
\label{eq:Letap,Linfty_estimate_thm1}
\begin{aligned}
	\|\partial_t^\ell u^\lin(t,\cdot)\|_{L^{q}}&\leq C (1+t)^{-\frac{1}{p_c}}\|f\|_{L^1\cap L^r};
\end{aligned}
\end{equation}
 assume that the function $\mu=\mu(s)$ fulfills hypothesis (i), (ii) and (iv) together with the integral condition \eqref{eq:muint}.
Thus, there exits $\eps_0\in (0,\epsb)$ such that if $f\in L^1\cap L^r$ satisfies
\begin{equation*} 
\label{eq:data_smallness_thm1}
\|f\|_{L^1\cap L^r}\leq \eps_0,
\end{equation*}
then,  there exists a unique solution $u$ to \eqref{eq:CPnon-linear} in  $ W_\loc^{\ell,\infty}([0,\infty), L^{p_c}\cap L^\infty)$.
\end{thm}
\begin{rem}
According to the properties of the differential operator $L$ one can be more precise about the functional space to which the solution $u$ to \eqref{eq:CPnon-linear} belongs, giving additionally information about the regularity of $u$ (in time and space) in appropriate Sobolev spaces; however, this is out of the scope of this paper.
\end{rem}
\begin{rem}
\label{rem:lower_derivative_estimate}
Suppose that the solution $u^\lin=u^\lin(t,x)$ to \eqref{eq:CPlinear} satisfies estimate \eqref{eq:Letap,Linfty_estimate_thm1} for some natural number $\ell $ in $\{1,\dots, m-1\}$; then, for any $k=0,1,\dots \ell-1$ we may estimate 
\begin{equation*} 
\label{eq:estimate_lower_derivatives1}
\|\partial_t^k u^\lin(t,\cdot)\|_{L^q}\leq C (1+t)^{-\frac{1}{p_c}+\ell-k}\|f\|_{L^1\cap L^r},
\end{equation*}
for some constant $C>0$ independent of $f$ and $t$; indeed, if  $k=\ell-1$ then, by applying the Minkowsky's integral inequality, we may estimate 
\begin{align*}
\label{eq:estimate_lower_derivatives2}
\|\partial_t^{k}u^\lin(t,\cdot)\|_{L^q}& \lesssim \int_0^t \|\partial_t^{\ell}u^\lin(\tau,\cdot)\|_{L^q}\,d\tau+\|\partial_t^{\ell-1}u^\lin(0,\cdot)\|_{L^q}\\ & \lesssim  (1+t)^{-\frac{1}{p_c}+1}\|f\|_{L^1\cap L^r},
\end{align*}
being $\partial_t^{\ell-1}u^\lin(0,\cdot)=0$.
Similarly, if $k\leq \ell-2$ we may estimate
\begin{align*}
\|\partial_t^k u(t,\cdot)\|_{L^q}&\lesssim \int_0^t \int_0^{t_1} \dots \int_0^{t_{\ell-k-1}} \|\partial_t^{\ell}u(\tau,\cdot)\|_{L^q}\,d\tau\,dt_{\ell-k-1}\,\dots dt_1\\&+\sum_{j=0}^{\ell-1-k}t^{\ell-1-k-j}\|\partial_t^{\ell-1-j}u(0,\cdot)\|_{L^q
}%
\lesssim (1+t)^{-\frac{1}{p_c}+\ell-k}\|f\|_{L^1\cap L^r},
\end{align*}
being $\partial_t^j u^\lin(0,x)=0$ for all $j=0,\dots, \ell-1$.\\
\end{rem}
\begin{proofof}{Theorem \ref{thm:existence}}
Let $R\leq \epsb/2$ with $\epsb>0$ given in hypothesis (ii) and (iv). We introduce the solution space 
\[ X= \{ u\in W^{\ell,\infty}_{\loc}([0,\infty), L^{p_c} \cap L^\infty) : \|u\|_{X}\leq R\},\]
where 
\[ \|u\|_{X}:= \sup_{s\in [0,\infty)}\big\{ (1+s)^{\frac{1}{p_c}}\|\partial_t^\ell u(s,\cdot)\|_{L^{p_c}\cap L^\infty}\big\}+\|u\|_{X_0},\]
where $\|\cdot\|_{X_0}=0$ if $\ell=0$, whereas 
\begin{equation*}
\|u\|_{X_0}:=\sup_{s\in [0,\infty)}\Big\{ \sum_{k=0}^{\ell-1}(1+s)^{\frac{1}{p_c}+k-\ell}\|\partial_t^k u(s,\cdot)\|_{L^{p_c}\cap L^\infty}\Big\},
\end{equation*}
if $\ell\geq 1$.
A function $u\in X$ is a solution to \eqref{eq:CPnon-linear} if, and only if, it satisfies identity \eqref{eq:duhamel}, i.e.
\begin{equation} 
\label{eq:integral_equation}
u(t,x)= u^\lin(t,x)+Nu(t,x)
\end{equation}
in $X$, where $N$ is the non-linear integral operator defined by
\begin{equation}
\label{eq:N}
Nu(t,x) : = \int_0^tK(t-s, \cdot)\ast_{(x)} F(\partial_t^\ell u(s,\cdot))\,ds.
\end{equation}
Since \eqref{eq:Letap,Linfty_estimate_thm1} holds true, we obtain
\begin{equation}
\label{eq:ulin_boundness}
\|u^\lin\|_{X} \leq C_1 \|f\|_{L^1\cap L^r},
\end{equation} 
for some $C_1>0$, independent of $f$ (see also Remark \ref{rem:lower_derivative_estimate}).\\
In the following, we will prove that there exists $C_2>0$ such that 
\begin{equation}
\label{eq:N_lipshitz}
\|Nu-Nv\|_{X}\leq C_2 \|u-v\|_X (\|u\|_{X}^{p_c-1}+\|v\|_{X}^{p_c-1}),
\end{equation}
for any $u,v\in X$.
Due to the definition of $\|\cdot\|_X$ we have,
\begin{equation*}
\|\partial_t^\ell u(s,\cdot)\|_{L^q}\leq (1+s)^{-\frac{1}{p_c}}\|u\|_{X}, \quad \text{ for any } q\in [p_c,\infty].
\end{equation*}
Thanks to estimate \eqref{eq:Letap,Linfty_estimate_thm1} we may estimate for any $q\in [p_c,\infty]$
\begin{equation}
\label{eq:Nu-Nv_estimate}
\begin{aligned}
\|\partial_t^\ell (Nu&(t,\cdot)-Nv(t,\cdot))\|_{L^q} \\ &\leq C \int_0^t (1+t-s)^{-\frac{1}{p_c}}\| F(\partial_t^\ell u(s,\cdot))-F(\partial_t^\ell v(s,\cdot))\|_{L^1\cap L^r}\,ds.
\end{aligned} 
\end{equation}
Moreover, we have
\[ |\partial_t^\ell u(s,x)|\leq \|\partial_t^\ell u(s,\cdot)\|_{L^\infty}\leq (1+s)^{-\frac{1}{p_c}}\|u\|_X\leq R,\]
for all $s\geq 0$. Therefore, since $R<\epsb$ and $F$ satisfies \eqref{eq:Fcontraction_intro}, applying H\"older inequality, we find 
\begin{align*}
\|F(\partial_t^\ell & u(s,\cdot))-F(\partial_t^\ell v(s,\cdot))\|_{L^\omega} \\ &\lesssim \|\mu((|\partial_t^\ell u|+|\partial_t^\ell v|)(s,\cdot))\|_{L^\infty} \\& \hspace{3em} \times \|\partial_t^\ell(u-v)(s,\cdot)(|\partial_t^\ell u|+|\partial_t^\ell v|)^{p_c-1}(s,\cdot)\|_{L^\omega} \\
& \lesssim \|\mu((|\partial_t^\ell u|+|\partial_t^\ell v|)(s,\cdot))\|_{L^\infty} \\& \hspace{3em} \times \|\partial_t^\ell(u-v)(s,\cdot)\|_{L^{\omega p_c}}(\|\partial_t^\ell u(s,\cdot)\|_{L^{\omega p_c}}^{p_c-1}+\|\partial_t^\ell v(s,\cdot)\|_{L^{\omega p_c}}^{p_c-1}) ,
\end{align*}
for any $\omega\geq 1$. In particular,  it holds
\begin{align*} 
\|\partial_t^\ell (u-v)(s,\cdot)\|_{L^{\omega p_c}}&(\|\partial_t^\ell u(s,\cdot)\|_{L^{\omega p_c}}^{p_c-1}+\|\partial_t^\ell v(s,\cdot)\|_{L^{\omega p_c}}^{p_c-1}) \\
\leq &(1+s)^{-1}\|\partial_t^\ell (u-v)\|_X(\|\partial_t^\ell u\|_X^{p_c-1}+\|\partial_t^\ell v\|_X^{p_c-1}).
\end{align*}
Moreover, since $\mu$ is non-decreasing we may estimate 
\begin{align*}
\|\mu(|\partial_t^\ell u(s,\cdot)|&+|\partial_t^\ell v(s,\cdot)|)\|_{L^\infty} \\ & \leq \mu\big((\|u\|_X+\|v\|_X)(1+s)^{-\frac{1}{p_c}}\big)\leq \mu\big(2R (1+s)^{-\frac{1}{p_c}}\big).
\end{align*}
We conclude that there exists $C>0$, independent of $f$ such that
\begin{align*}
\|F(\partial_t^\ell u(s,\cdot))& -F(\partial_t^\ell v(s,\cdot))\|_{L^\omega} \\ & \leq C(1+s)^{-1} \mu\big(2R (1+s)^{-\frac{1}{p_c}} )\|u-v\|_X(\|u\|_X^{p_c-1}-\|v\|_X^{p_c-1}).
\end{align*}
Then, we may estimate
\begin{align*}
\|\partial_t^\ell (Nu(t,\cdot)-&Nv(t,\cdot))\|_{L^q}\leq C \|u-v\|_X(\|u\|_X^{p_c-1}-\|v\|_X^{p_c-1})\\
& \times \int_0^t (1+t-s)^{-\frac{1}{p_c}}(1+s)^{-1} \mu\big(2R (1+s)^{-\frac{1}{p_c}} )\,ds.
\end{align*} 
We split the integral in $[0,t/2]$ and $[t/2,t]$. On the one hand,
\begin{align*}
\int_{\frac{t}{2}}^t \dots ds & \leq \Big(\sup_{\tau\in [0,\epsb]}\mu(\tau)\Big) (1+t/2)^{-1} \int_\frac{t}{2}^t (1+t-s)^{-\frac{1}{p_c}}\,ds \\
& \leq C \Big(\sup_{\tau\in [0,\epsb]}\mu(\tau)\Big) (1+t)^{-\frac{1}{p_c}};
\end{align*}
on the other hand, being $R\leq \epsb/2$ we have
\begin{align*}
\label{eq:int_0_t/2}
\int_0^{\frac{t}{2}} \dots ds \leq (1+t/2)^{-\frac{1}{p_c}} \int_0^{\frac{t}{2}} (1+s)^{-1} \mu\big(\epsb (1+s)^{-\frac{1}{p_c}} )\,ds\leq  C (1+t)^{-\frac{1}{p_c}};
\end{align*}
indeed, the integral term is uniformly bounded with respect to $t\in [0,\infty)$ thanks to the integral condition \eqref{eq:muint}.  Thus, we conclude
\begin{equation*}
\|\partial_t^\ell (Nu(t,\cdot)-Nv(t,\cdot))\|_{L^q}\leq C  (1+t)^{-\frac{1}{p_c}}\|u-v\|_X(\|u\|_X^{p_c-1}-\|v\|_X^{p_c-1}),
\end{equation*}
for any $q\in [p_c,\infty]$. \\
The proof of inequality \eqref{eq:N_lipshitz} directly follows if $\ell=0$; if $\ell\geq 1$ applying the same approach as in Remark \ref{rem:lower_derivative_estimate} we can easily derive for any $q\in [p_c,\infty]$ the estimate 
\begin{equation*}
\|\partial_t^k (Nu(t,\cdot)-Nv(t,\cdot))\|_{L^q}\leq C  (1+t)^{-\frac{1}{p_c}+\ell-k}\|u-v\|_X(\|u\|_X^{p_c-1}-\|v\|_X^{p_c-1}),
\end{equation*}
for every $k=0,\dots, \ell-1$. Then, there exists $C_2>0$ such that inequality \eqref{eq:N_lipshitz} holds true.

The proof of the global in time existence of small data solutions is obtained applying a classic contraction argument. We define
\begin{equation}
\label{eq:R_choice1}
R=2C_1\,\|f\|_{L^1\cap L^r}, 
\end{equation}
where $C_1$ is as in~\eqref{eq:ulin_boundness}. We consider $R\leq \epsb/2$, sufficiently small, such that 
\begin{equation}
\label{eq:R_choice2}
2C_2R^{p_c-1}\leq 1/2,
\end{equation}
where $C_2$ is as in~\eqref{eq:N_lipshitz}. By~\eqref{eq:ulin_boundness} and~\eqref{eq:N_lipshitz}, it follows that the operator $u^\lin(t,x) + N$ maps $X$ into itself; in particular, due to \eqref{eq:N_lipshitz} and \eqref{eq:R_choice2}, it is a contraction.

For any arbitrarily large~$T>0$ we can use the same approach to prove that inequalities \eqref{eq:ulin_boundness} and \eqref{eq:N_lipshitz} holds true, for some constant $C_1$ and $C_2$ independent of $T$, replacing $X$ by the Banach space $X(T)$, defined by
\[ X(T)= \{ u\in W^{\ell,\infty}([0,T], L^{p_c} \cap L^\infty) : \|u\|_{X(T)}\leq R\},\]
where 
\[ \|u\|_{X(T)}:= \sup_{s\in [0,T]}\big\{ (1+s)^{\frac{1}{p_c}}\|\partial_t^\ell u(s,\cdot)\|_{L^{p_c}\cap L^\infty}\big\}+\|u\|_{X_0(T)},\]
with $\|\cdot\|_{X_0(T)}=0$ if $\ell=0$, whereas 
\begin{equation*}
\|u\|_{X_0(T)}:=\sup_{s\in [0,T]}\Big\{ \sum_{k=0}^{\ell-1}(1+s)^{\frac{1}{p_c}+k-\ell}\|\partial_t^k u(s,\cdot)\|_{L^{p_c}\cap L^\infty}\Big\},
\end{equation*}
if $\ell\geq 1$. Since $X(T)$ is a Banach space, this allows to conclude that  there is unique fixed point for the operator $u^\lin(t,x) + N$ in $X(T)$, that is, a unique solution to~\eqref{eq:integral_equation} in $X(T)$. Being $T$ arbitrary large, we get that a unique solution exists in $X$. 
\end{proofof}
\subsection{The case $p_c=1$}
\label{sec:existence-p=1}
The proof of Theorem \ref{thm:existence} heavily relies on the contractivity of the integral operator $N$ defined by \eqref{eq:N}; in the case $p_c>1$ such property follows directly by inequality \eqref{eq:N_lipshitz}, taking $\|u\|_X$, $\|v\|_X\leq R$, with $R$ defined by \eqref{eq:R_choice1} sufficiently small to satisfy inequality \eqref{eq:R_choice2}. Clearly, the same idea can not be used if $p_c=1$; in order to apply a contraction argument in this latter case we require that there exists $\bar{\eps}>0$ such that for all $\eps \in [0,\bar{\eps}]$ the non-linear term $F$ satisfies 
\begin{equation} 
\label{eq:fmu_p=1}
F(0)=0, \quad |F(u)-F(v)|\leq C_\eps |u-v| \quad \text{for all }|u|,|v|\in [0,\eps],
\end{equation}
for some $C_\eps>0$ that tends to $0$ as $\eps\to 0$.
\begin{ex}
\label{example:F_pc=1}
For instance, one can consider $F:\R\to \R$ such that
\[ F(u)=|u|\mu(|u|),\]
where $\mu$ is a continuous function such that $\mu(0)=0$ and satisfying properties (i), (ii) and (iv) with $p_c=1$. \\
Indeed, thanks to assumption (iv) on $\mu$, there exists $\bar{\eps}>0$ such that the function  $F(s)=|s|\mu(|s|)$ satisfies
\[F(0)=0, \quad |F(u)-F(v)|\leq C|u-v|\mu(|u|+|v|), \]
for $|u|<\epsb$ and $ |v|<\epsb$, for some constant $C>0$ independent of $u$ and $v$. In particular, since $\mu$ is non-decreasing, for every $\eps\in [0,\epsb]$ 
it holds
\[  |F(u)-F(v)|\leq C|u-v|\mu(2\eps),\]
for all $|u|, |v|\in [0,\eps]$; namely, for every $\eps\in [0,\epsb]$  and $|u|, |v|\in [0,\eps]$ the desired inequality \eqref{eq:fmu_p=1} holds true taking $C_\eps:=C\mu(2\eps)$; being $\mu=\mu(s)$ a continuous fucntion and $\mu(0)=0$ then $C_\eps\to 0$ as $\eps\to 0$.
\end{ex} 
\begin{prop}
\label{thm:existence2}
Let $\delta \in [1,\infty]$; suppose that there exists $r\in [\delta, \infty]$ such that for every $q\in [\delta ,\infty]$ and $f\in L^\delta \cap L^r$ the solution to the linear Cauchy problem \eqref{eq:CPlinear} satisfies 
\begin{equation}
	\label{eq:Letap,Linfty_estimate_thm2}
	\begin{aligned}
		\|\partial_t^\ell u^\lin(t,\cdot)\|_{L^{q}}&\leq C (1+t)^{-\beta}\|f\|_{L^\delta\cap L^r},
		%
	\end{aligned}
\end{equation}
for some $\beta > 1$. Also assume that the non-linearity $F$ verifies assumption \eqref{eq:fmu_p=1} for some $\bar{\eps}>0$ arbitrarily small.
Thus, there exists $\eps_0\in (0,\bar \eps)$ such that if $f\in L^\delta \cap L^r$ verifies
\begin{equation*} 
	\label{eq:data_smallness_thm2}
	\|f\|_{L^\delta\cap L^r}\leq \eps_0,
\end{equation*}
then, there exists a unique solution $u$ to \eqref{eq:CPnon-linear} in $ W^{\ell,\infty}_\loc([0,\infty), L^{\delta}\cap L^\infty)$.
\end{prop}
\begin{rem}
\label{rem:lower_derivative_estimate_p=1}
As already discussed in Remark \ref{rem:lower_derivative_estimate}, if the solution $u^\lin=u^\lin(t,x)$ to \eqref{eq:CPlinear} satisfies estimate \eqref{eq:Letap,Linfty_estimate_thm2} for some natural number $\ell$ in $\{1,\dots, m-1\}$, then, for any $k=0,1,\dots \ell-1$ and $q\in [\delta, \infty]$ we may estimate 
\begin{equation} 
\label{eq:estimate_lower_derivatives3}
\|\partial_t^k u^\lin(t,\cdot)\|_{L^q}\leq C (1+t)^{-\beta+\ell-k}\|f\|_{L^\delta\cap L^r},
\end{equation}
for some constant $C>0$ independent of $f$ and $t$.
\end{rem}
\begin{proofof}{Porposition \ref{thm:existence2}}
We proceed similarly as in the proof of Theorem \ref{thm:existence}; we introduce the solution space 
\[ X= \{ u\in W^{\ell,\infty}_{\loc}([0,\infty), L^{\delta} \cap L^\infty) : \|u\|_{X}\leq R\},\]
where $0<R<\bar \eps$, with $\bar \eps$ given in \eqref{eq:fmu_p=1} and
\begin{align*}
	\|u\|_{X}:= \sup_{s\in [0,\infty)} 
	\Big\{(1+s)^{\beta}\|\partial_t^\ell u(s,\cdot)\|_{L^{\delta}\cap L^\infty}\Big\}+\|u\|_{X_0},
\end{align*}
where $\|\cdot\|_{X_0}=0$ if $\ell=0$, whereas 
\begin{equation*}
\|u\|_{X_0}:=\sup_{s\in [0,\infty)}\Big\{ \sum_{k=0}^{\ell-1}(1+s)^{\beta+k-\ell}\|\partial_t^k u(s,\cdot)\|_{L^{\delta}\cap L^\infty}\Big\},
\end{equation*}
if $\ell\geq 1$. A function $u\in X$ is a solution to \eqref{eq:CPnon-linear} if, and only if, it satisfies \eqref{eq:duhamel} in $X$, i.e.
\begin{equation} 
	\label{eq:integral_equation1}
	u(t,x)=u^\lin(t,x)+Nu(t,x),
\end{equation}
in $X$, where $N$ is the non-linear integral operator defined by
\[ Nu(t,x) : = \int_0^t K(t-s,\cdot)\ast_{(x)} F(\partial_t^\ell u(s,\cdot))\,ds.\]

Since $f\in L^\delta\cap L^r$ and \eqref{eq:Letap,Linfty_estimate_thm2} holds true, we may estimate
\begin{equation}
	\label{eq:ulin_boundness1}
	\|u^\lin\|_{X} \leq C_1 \|f\|_{L^\delta\cap L^r},
\end{equation} 
for some $C_1>0$, independent of $f$ (see also Remark \ref{rem:lower_derivative_estimate_p=1}). 
In the following,  we will prove that there exists $C'_R>0$ such that 
\begin{equation}
	\label{eq:N_lipshitz_p=1}
	\|Nu-Nv\|_{X}\leq C'_R \|u-v\|_X,
\end{equation}
for some constant $C'_R>0$ that tends to $0$ as $R\to 0$.\\
Due to the definition of $\|\cdot\|_X$ we have,
\begin{equation*}
	\|\partial_t^\ell u(s,\cdot)\|_{L^q}\leq (1+s)^{-\beta }\|u\|_{X}, \quad \text{ for any } q\in [\delta ,\infty].
\end{equation*}
%
%
Morover, thanks to estimate \eqref{eq:Letap,Linfty_estimate_thm2} we may estimate for any $q\in [\delta ,\infty]$
\begin{align*}
	\|\partial_t^\ell (Nu(t,\cdot)&-Nv(t,\cdot))\|_{L^q}\\&\leq C \int_0^t (1+t-s)^{-\beta}\|F(\partial_t^\ell u(s,\cdot))-F(\partial_t^\ell v(s,\cdot))\|_{L^\delta\cap L^r}\,ds.
\end{align*} 
Moreover, we have
\[ |\partial_t^\ell u(s,x)|\leq \|\partial_t^\ell u(s,\cdot)\|_{L^\infty}\leq (1+s)^{-\beta}\|u\|_X\leq R,\]
for all $s\geq 0$. Therefore, since $R<\bar\eps$ and $F$ satisfies \eqref{eq:fmu_p=1}, applying H\"older inequality, we find 
\begin{align*}
	\|F(\partial_t^\ell u(s,\cdot))-F(\partial_t^\ell v(s,\cdot))\|_{L^\omega} \leq C_R \|\partial_t^\ell (u-v)(s,\cdot)\|_{L^{\omega }},
\end{align*}
for any $\omega\geq \delta$, where $C_R$ tends to $0$ as $R\to 0$. In particular, it holds
\[ \|\partial_t^\ell (u-v)(s,\cdot)\|_{L^{\omega }} \leq (1+s)^{-\beta}\|u-v\|_X.\]
We conclude 
\begin{align*}
	\|F(\partial_t^\ell u(s,\cdot))-F(\partial_t^\ell v(s,\cdot))\|_{L^\omega}\leq C_R(1+s)^{-\beta}  \|u-v\|_X,
\end{align*}
where $C_R$ tends to $0$ as $R\to 0$; then, since $\beta>1$, for any $q\geq \delta $ we may estimate
\begin{align*}
	\|\partial_t^\ell (Nu(t,\cdot)-Nv(t,\cdot))\|_{L^q}&\leq CC_R \|u-v\|_X\int_0^t (1+t-s)^{-\beta}(1+s)^{-\beta}\,ds \\
& \leq C'_R (1+s)^{-\beta} \|u-v\|_X,
\end{align*} 
for some $C'_R>0$ that tends to $0$ as $R$ tends to $0$.
The proof of inequality \eqref{eq:N_lipshitz_p=1} directly follows. 

Also in this case, the proof of the global in time existence of small data solutions is obtained applying a classic contraction argument: we define
\[ R=2C_1\,\|f\|_{L^\delta\cap L^r}, \]
where $C_1$ is as in~\eqref{eq:ulin_boundness1}. Since \eqref{eq:N_lipshitz_p=1} holds true and $C'_R\to 0$ as $R\to 0$, we can consider $R<\bar \eps$ sufficiently small, such that $C'_R \leq 1/2$.\\
By~\eqref{eq:ulin_boundness1} and~\eqref{eq:N_lipshitz_p=1}, it follows that the operator $K(t,\cdot)\ast_{(x)}f + N$ maps $X$ into itself; in particular, due to~\eqref{eq:N_lipshitz_p=1} and our choice of $R$, it is a contraction.

For any arbitrarily large~$T>0$, we can use the same approach to prove that inequalities \eqref{eq:ulin_boundness1} and \eqref{eq:N_lipshitz_p=1} holds true if we replace $X$ by the Banach space $X(T)$, defined by
\[ X(T) = \{ u\in W^{\ell, \infty}([0,T], L^{\delta }\cap L^\infty):\, \|u\|_{X(T)}\leq R\}, \]
where 
\begin{align*}
	\|u\|_{X(T)}:= \sup_{s\in [0,T]} 
	\Big\{(1+s)^{\beta}\|\partial_t^\ell u(s,\cdot)\|_{L^{\delta}\cap L^\infty}\Big\}+\|u\|_{X_0(T)},
\end{align*}
where $\|\cdot\|_{X_0}=0$ if $\ell=0$, whereas 
\begin{equation*}
\|u\|_{X_0(T)}:=\sup_{s\in [0,T]}\Big\{ \sum_{k=0}^{\ell-1}(1+s)^{\beta+k-\ell}\|\partial_t^k u(s,\cdot)\|_{L^{\delta}\cap L^\infty}\Big\},
\end{equation*}
if $\ell\geq 1$.
Since $X(T)$ is a Banach space, this allow to conclude that  there is a unique fixed point for $K(t,\cdot)\ast_{(x)}f + N$ in $X(T)$, that is, a unique solution to~\eqref{eq:integral_equation1} in $X(T)$. Being $T$ arbitrarily large, we get the existence of a unique solution in $X$. 
\end{proofof}
\section{Examples of application}\label{sec:examples}
In this section we provide examples of evolution models to which our main results can be applied. The equations which we take into account have been deeply studied in the literature and are obtained as a perturbation of well known equations of mathematical physics.
\subsection{The classical damped wave equation in space dimension $n\leq 3$.}
We provide the following result for a non-linear classical damped wave equation.
\begin{prop}
Consider for $n\geq 1$ the Cauchy problem 
\begin{equation}
\label{eq:dampedwave}
\begin{cases}
\partial_{t}^2 u- \Delta u + \partial_tu =|u|^{p_c}\mu(|u|), \quad (t,x)\in [0,+\infty)\times \R^n,\\
u(0,x)=0, \quad \partial_tu(0,x)=f(x),
\end{cases}
\end{equation}
where $p_c=1+2/n$ and $\mu:[0,\infty) \to [0,\infty)$ satisfies assumptions (i)-(iv). \\ Assume that $f\in L^1(\R^n)$ satisfies the sign condition 
\begin{equation*}
\label{eq:sign_f_dampedwave}
 \int_{\R^n} f(x)\,dx>0;
\end{equation*}
if there exists a global weak solution $u\in L^{p_c}_{\text{loc}}\big([0,\infty),L^{p_c}_{\text{loc}}(\R^n)\big)$ to \eqref{eq:dampedwave}, according to Definition \ref{def:weaksolution}, then
\begin{equation}
\label{eq:muint_dampedwave}
\int_{0}^{c_0} \frac{\mu(s)}{s}\,ds<\infty,
\end{equation}
for some $c_0>0$ sufficiently small. 

On the other hand, in space dimension $n\leq 3$, if $\mu=\mu(s)$ verifies the integral condition \eqref{eq:muint_dampedwave}, thus there exists $\eps_0>0$ such that for all $f\in L^1\cap L^2$ satisfying
\[\|f\|_{L^1\cap L^2}\leq \eps_0,\]
Cauchy problem \eqref{eq:dampedwave} admits a unique global (in time) solution $u$ in the space $L^\infty_{loc}([0,\infty), L^{p_c}\cap L^\infty).$
\end{prop}
The proof of the non-existence result is an immediate consequence of Theorem \ref{thm:non-existenceCauchy}.
In \cite{EGR} it was already proved its optimality in space dimension $n=1$; more in detail, if $\mu$ fulfils \eqref{eq:muint_dampedwave} then a unique global solution to \eqref{eq:dampedwave} exists in $C([0,\infty), L^2\cap L^\infty)$ assuming the initial data to be sufficiently small in $L^1\cap L^2$. The same is true in space dimension $n=2$ if one assumes additional regularity $H^1$ for the initial data. \\
The application of Theorem \ref{thm:existence} allows to reach an analogous conclusion also in space dimension $n=3$; indeed, we know by \cite{Nishihara2003} that the solution $u^\lin(t,\cdot)=u^\lin(t,\cdot)$ to the linear Cauchy problem 
\begin{equation}
	\label{eq:dampedwave_lin}
	\begin{cases}
		\partial_{t}^2 u- \Delta u + \partial_tu = 0, \\
		u(0,x)=0, \quad \partial_t u(0,x)=f(x),
	\end{cases}
\end{equation} 
satisfies for every $f\in L^1\cap L^\infty$ the long time decay estimate
%
%
%
%
%
%
%
\[ \|u^\lin(t,\cdot)\|_{L^q} \lesssim (1+t)^{-\frac{3}{2}(1-\frac{1}{q})}\|f\|_{L^1\cap L^\infty},\]
for any $q\in [1,\infty]$; in particular, for all $q\geq 5/3$ it holds 
\[ \|u^\lin(t,\cdot)\|_{L^q}\lesssim (1+t)^{-\frac{3}{5}}\|f\|_{L^1\cap L^\infty}. \]
 Thus, if $\mu$ satisfies \eqref{eq:muint_dampedwave} for some $c_0>0$ then, Theorem \ref{thm:existence} guarantees that Cauchy problem \eqref{eq:dampedwave} admits a unique global solution $u\in L^\infty_\loc([0,\infty), L^\frac{5}{3}\cap L^\infty)$, provided that $f\in L^1\cap L^\infty$ is sufficiently small.
\subsection{The structurally damped $\sigma$-evolution equation}
\label{sec:structurally}
We provide the following result for a non-linear damped $\sigma$-evolution equation.
\begin{prop}
Consider for $n\geq 1$ the Cauchy problem 
\begin{equation}
\label{eq:sigma-evolution_ell=0}
\begin{cases}
\partial_t^2 u + (-\Delta)^\delta \partial_t u + (-\Delta)^\sigma u = | u|^{p_c}\mu( |u|),  \quad (t,x)\in [0,\infty)\times \R^n,  \\
u(0,x)=0, \quad \partial_t u(0,x)=f(x),&
\end{cases}
\end{equation}
where $0\leq \delta< \sigma$ are natural numbers, $\mu:[0,\infty) \to [0,\infty)$ satisfies assumptions (i)-(iv), and $p_c$ is determined by formula \eqref{eq:pcritical}, i.e. 
\begin{equation*}
\begin{aligned}
 p_c&=1+\frac{2\sigma}{n}, &&\quad \text{ if } \delta=0 &&&\textit{(classical damping)},\\
p_c&=1+\frac{2\sigma}{n-2\delta}, &&\quad \text{ if } 2\delta\in (0,\sigma) &&& \textit{(structural effective damping)},
\end{aligned}
\end{equation*}
provided that $n>2\delta$, and 
\[ p_c=1+\frac{2\sigma}{n-\sigma}, \hspace{1.7em}\quad \text{ if } 2\delta\geq \sigma\hspace{2.3em}  \textit{(structural non-effective damping)},\]
provided that $n>\sigma$. \\ Assume that $f\in L^1(\R^n)$ satisfies the sign condition 
\begin{equation*}
\label{eq:sign_f_sigmaevol}
 \int_{\R^n} f(x)\,dx>0;
\end{equation*}
if there exists a global weak solution $u\in L^{p_c}_{\text{loc}}\big([0,\infty),L^{p_c}_{\text{loc}}(\R^n)\big)$ to \eqref{eq:sigma-evolution_ell=0}, according to Definition \ref{def:weaksolution}, then
\begin{equation}
\label{eq:muint_sigmaevol}
\int_{0}^{c_0} \frac{\mu(s)}{s}\,ds<\infty,
\end{equation}
for some $c_0>0$, sufficiently small. 

On the other hand, suppose that $2\delta \in [0,\sigma)$ and $n\in (2\delta, 2\sigma)$, or $2\delta \in [\sigma,2\sigma]$ with $\sigma>1$ and $n\in (\sigma,2\sigma)$; if $\mu=\mu(s)$ verifies the integral condition \eqref{eq:muint_sigmaevol}, then there exists $\eps_0>0$ such that for all $f\in L^1\cap L^2$ satisfying
\[\|f\|_{L^1\cap L^2}\leq \eps_0,\]
Cauchy problem \eqref{eq:sigma-evolution_ell=0} admits a unique global (in time) solution $u$ in the space $L^\infty_{loc}([0,\infty), L^{p_c}\cap L^\infty).$
\end{prop}
The proof of the non-existence result is an immediate consequence of Theorem \ref{thm:non-existenceCauchy} (see also Example \ref{example:pc}). The existence of global small data solutions for the same problem with power non-linearity has been deeply investigated in the literature (see \cite{DAE2014nonlinear, DAE2017, DAE22, DAR2014}) with the application of suitable $L^p-L^q$ long time decay estimates, whose structure is related to the different influence of the damping term $(-\Delta)^{\delta} \partial_tu$, when it is \textit{classical}, \textit{effective} or \textit{non-effective}  (see also the classification of effective/non-effective damping introduced in \cite{DAE2016}). If the damping is classical or effective (namely, $2\delta\in [0,\sigma)$) and $n<2\sigma$, employing the estimates obtained in \cite{DAE2017} one can prove that the solution $u^\lin=u^\lin(t,x)$ to the linear Cauchy problem 
\begin{equation}
\label{eq:sigma-evolution_ell=0_lin}
\begin{cases}
\partial_t^2 u + (-\Delta)^\delta \partial_t u + (-\Delta)^\sigma u = 0& \\
u(0,x)=0, \quad \partial_tu(0,x)=f(x),&
\end{cases}
\end{equation}
satisfies
\[ \|u^\lin(t,\cdot)\|_{L^q}\lesssim (1+t)^{-\frac{1}{2(\sigma-\delta)}(n(1-\frac{1}{q})-2\delta)}\|f\|_{L^1\cap L^2},\]
for any $q\in [1,\infty]$. In particular, if $n\in (2\delta, 2\sigma)$ and  $q>p_c=1+2\sigma/(n-2\delta)$ then 
\[ \|u^\lin(t,\cdot)\|_{L^q}\lesssim (1+t)^{-\frac{1}{p_c}}\|f\|_{L^1\cap L^2},\]
for any $q\in [1,\infty]$. Thus, in this case Theorem \ref{thm:existence} guarantees the existence of a unique global small data solution to \eqref{eq:sigma-evolution_ell=0} in $L^\infty_\loc([0,\infty), L^{p_c}\cap L^\infty)$ if $\mu$ satisfies the integral condition \eqref{eq:muint_sigmaevol} for some $c_0>0$. \\
Similarly, in the non-effective damped case (namely, $2\delta \in [\sigma, 2\sigma]$) if $\sigma>1$ and $n\in (\sigma, 2\sigma)$, then a global solution to \eqref{eq:sigma-evolution_ell=0} exists if $\mu$ fulfils \eqref{eq:muint_sigmaevol} and  $f$ is small in $L^1\cap L^2$; this latter result, was already proved in \cite{DAG2024}. In particular, in the homogeneous case ($2\delta=\sigma$) the same result can be obtained in any space dimension, assuming small initial data in $L^1\cap L^{\frac{n}{\sigma}}$. 
Our Theorem \ref{thm:non-existenceCauchy} allows to conclude that the results obtained in \cite{DAG2024} are optimal.\\
As discussed in \cite{DAG2024}, in the case $\sigma=1$ and $2\delta \in (1,2]$ the wave structure appears and oscillations come into play; as a consequence, only a partial result about the critical power can be proved, for instance in \cite{DAR2014}, but the obtained existence exponent is far from being optimal. Recently, in \cite{DAE23} the authors provided sharp $L^p-L^q$ estimates for the corresponding linear model which can be applied to improve the known results about global (in time) existence.\\
Notice that the critical exponent $p_c$ for problem \eqref{eq:sigma-evolution_ell=0} remains the same if one consider fractional exponents $\delta$ and $\sigma$; this can be proved employing a different test function which allows to overcome the difficulties due to the non-local behavior of fractional differential operators (see \cite{DAF2021}).
\subsection{The structurally damped $\sigma$-evolution equation with derivative-type non-linearity}
We provide the following result for a non-linear damped $\sigma$-evolution equation with derivative-type non-linearity.
\begin{prop}
Consider for $n\geq 1$ the Cauchy problem 
\begin{equation}
\label{eq:sigma-evolution_ell=1}
\begin{cases}
\partial_t^2 u + (-\Delta)^\delta \partial_t u + (-\Delta)^\sigma u = |\partial_t u|^{p_c}\mu( |\partial_t u|), &\\
u(0,x)=0, \quad \partial_t u(0,x)=f(x),&
\end{cases}
\end{equation}
for $(t,x)\in [0,\infty)\times \R^n $, where $0< \delta< \sigma$ are natural numbers, $\mu:[0,\infty) \to [0,\infty)$ satisfies assumptions (i)-(iv), and $p_c$ is determined by formula \eqref{eq:pcritical}, i.e. 
\begin{equation*}
\begin{aligned}
%
p_c&=1+\frac{2\delta}{n}, &&\quad \text{ if } 2\delta\in (0,\sigma) &&& \textit{(structural effective damping)},
\end{aligned}
\end{equation*}
and 
\[ p_c=1+\frac{\sigma}{n}, \hspace{1.7em}\quad \text{ if } 2\delta\geq \sigma\hspace{2.3em}  \textit{(structural non-effective damping)}.\]
Assume that $f\in L^1(\R^n)$ satisfies the sign condition 
\begin{equation*}
\label{eq:sign_f_sigmaevol_ell=1}
 \int_{\R^n} f(x)\,dx>0;
\end{equation*}
if there exists a global weak solution $u\in W^{1,p_c}_{\text{loc}}\big([0,\infty),L^{p_c}_{\text{loc}}(\R^n)\big)$ to \eqref{eq:sigma-evolution_ell=1}, according to Definition \ref{def:weaksolution}, then
\begin{equation}
\label{eq:muint_sigmaevol_ell=1}
\int_{0}^{c_0} \frac{\mu(s)}{s}\,ds<\infty,
\end{equation}
for some $c_0>0$ sufficiently small. 

On the other hand, suppose that $2\delta=\sigma$, or $2\delta \in (\sigma,2\sigma)$ for some $\sigma\geq 3$ and $n< \sigma-2$; if $\mu=\mu(s)$ verifies the integral condition \eqref{eq:muint_sigmaevol_ell=1}, then there exists $\eps_0>0$ such that for all $f\in L^1\cap L^\infty$ satisfying
\[\|f\|_{L^1\cap L^\infty}\leq \eps_0,\]
Cauchy problem \eqref{eq:sigma-evolution_ell=1} admits a unique global (in time) solution $u$ in the space $W^{1,\infty}_{loc}([0,\infty), L^{1+\frac{\sigma}{n}}\cap L^\infty).$
\end{prop}
The proof of the non-existence result is an immediate consequence of Theorem \ref{thm:non-existenceCauchy} (see also Example 6.4 in \cite{DAF2021}). In \cite{DAE2017} and \cite{DAE21} the authors provided long-time decay estimates for the solution $u^\lin=u^\lin(t,x)$ to the corresponding linear problem 
\begin{equation}
\label{eq:sigma-evolution_ell=1_linear}
\begin{cases}
\partial_t^2 u + (-\Delta)^\delta \partial_t u + (-\Delta)^\sigma u = 0 \\
u(0,x)=0, \quad \partial_t u(0,x)=f(x);&
\end{cases}
\end{equation}
in the non-effective case $2\delta\geq \sigma$, for all $q\in [p_c,\infty]$ the solution to \eqref{eq:sigma-evolution_ell=1_linear} satisfies 
\[ \|\partial_t u^\lin(t,\cdot)\|_{L^q}\lesssim (1+t)^{-\frac{n}{\sigma}(1-\frac{1}{q})}\|f\|_{L^1\cap L^q},\]
provided that $n<\sigma-2$ if $2\delta>\sigma$. In particular, for all $q\geq p_c=1+\sigma/n$ the solution to \eqref{eq:sigma-evolution_ell=1_linear} satisfies 
\[ \|\partial_t u^\lin(t,\cdot)\|_{L^q}\lesssim (1+t)^{-\frac{1}{p_c}}\|f\|_{L^1\cap L^\infty}. \] 
In this case, as a consequence of Theorem \ref{thm:existence} a unique global solution to \eqref{eq:sigma-evolution_ell=1} exists in $W^{1,\infty}_\loc([0,\infty), L^{1+\frac{\sigma}{n}}\cap L^\infty)$ provided that $f$ is sufficiently small in $L^1\cap L^\infty$ and $\mu$ satisfies the integral condition \eqref{eq:muint_sigmaevol_ell=1}. 
\subsection{Higher order dissipative hyperbolic equations}
Let us consider a partial differential operator $Q(\partial_t,\partial_x)$ which is the sum of two or three homogeneous hyperbolic operators of order $m$ and $m-1$, or, respectively, $m$, $m-1$ and $m-2$; namely, $Q=Q_1$ or $Q=Q_2$ where
\begin{equation}
\label{eq:Q_case1}
Q_1(\partial_t,\partial_x)=\PI_m(\partial_t,\partial_x)+\PI_{m-1}(\partial_t,\partial_x),
\end{equation} 
or 
\begin{equation}
\label{eq:Q_case2} Q_2(\partial_t,\partial_x)=\PI_m(\partial_t,\partial_x)+\PI_{m-1}(\partial_t,\partial_x)+\PI_{m-2}(\partial_t,\partial_x),
\end{equation}
with $ \PI_k(\partial_t,\partial_x)$ homogenous hyperbolic operators of order $k$, defined by
\[ \PI_k(\partial_t,\partial_x):= \sum_{j+|\alpha|=k} c_{j,\alpha} \partial_t^j \partial_x^\alpha, \quad k=m,\, m-1\, m-2,\]
with corresponding symbol 
\begin{equation*}
\PI_k(\lambda,\xi):= \sum_{j+|\alpha|=k} c_{j,\alpha}\lambda^j \xi^\alpha, \quad k=m-2, \, m-1,\, m.
\end{equation*}
In particular, we fix $\varrho=0$ if $Q=Q_1$ or $\varrho=1$ if $Q=Q_2$. Moreover, without loss of generality we suppose $c_{m,0}=1$ and $c_{m-1,0}>0$ if $Q=Q_1$  and, additionaly, $c_{m-2,0}>0$ if $Q=Q_2$.

The critical exponent for the semi-linear problem
\[ Q(\partial_t,\partial_x)=|u|^p,\] is strictly related to the behavior of the roots of the polynomials $\PI_m(\lambda,\xi)$, $\PI_{m-1}(\lambda,\xi)$ and $\PI_{m-2}(\lambda,\xi)$. Let us recall the following definitions.
\begin{defn}
A complex polynomial $p(z)$ is \textit{hyperbolic} if its roots are real numbers. It is\textit{ strictly hyperbolic} if its roots are real and simple. 
\end{defn}
\begin{defn}
\label{def:strictly_interlacing}
Let $p_{k-1}(z)$ be an hyperbolic polynomial of degree $k-1$ with roots $\tilde{\lambda}_1\leq \tilde{\lambda}_2\leq \dots \leq \tilde{\lambda}_{m-1}$, and let $p_m(z)$ be an hyperbolic polynomial of degree $m$ with roots $\lambda_1\leq \dots \leq \lambda_m$.\\
We say that $p_{k-1}(z)$ and $p_k(z)$\textit{ interlace} if
\[ \lambda_1 \leq \tilde \lambda_1 \leq \lambda_2\leq \tilde{\lambda}_2 \leq \dots \leq \tilde{\lambda}_{m-1}\leq \lambda_m.\]
In particular, we say that $p_{k-1}(z)$ and $p_k(z)$ \textit{strictly interlace} if
\[ \lambda_1 <\tilde \lambda_1 < \lambda_2< \tilde{\lambda}_2 < \dots < \tilde{\lambda}_{m-1}< \lambda_m.\]
\end{defn}
In this section we assume that the following conditions are satisfied by $\PI_m(\lambda,\xi)$ and $\PI_{m-1}(\lambda, \xi)$, if $Q=Q_1$.  
\begin{hypothesis} \label{hyp:1}For all $\xi'\in S^{n-1}:=\{\xi'\in \R^{n-1}: |\xi'|=1\}$:
\begin{itemize}
\item $\PI_{m}(\lambda, \xi')$ and $\PI_{m-1}(\lambda, \xi')$ are strictly hyperbolic;
\item $\PI_{m-1}(\lambda, \xi')$ and $\PI_{m}(\lambda, \xi')$ strictly interlace;
\item $\PI_{m-1}(0, \xi')=0$. 
\end{itemize}
\end{hypothesis}
\noindent Similarly, if $Q=Q_2$ we suppose that the following hypothesis holds true:
\begin{hypothesis}\label{hyp:2} For all $\xi'\in S^{n-1}$:
\begin{itemize}
\item $\PI_{m}(\lambda, \xi')$, $\PI_{m-1}(\lambda, \xi')$ and $\PI_{m-2}(\lambda, \xi')$ are strictly hyperbolic;
\item $\PI_{m-1}(\lambda, \xi')$ and $\PI_{m}(\lambda, \xi')$ strictly interlace;
\item $\PI_{m-2}(\lambda, \xi')$ and $\PI_{m-1}(\lambda, \xi')$ strictly interlace;
\item $\PI_{m-2}(0, \xi')=0$. 
\end{itemize}
\end{hypothesis}
\noindent The application of Theorems \ref{thm:non-existenceCauchy} and  \ref{thm:existence} allows to prove the following result.

\begin{prop}
\label{prop:higherorder}
Let $m\geq 3$ a natural number and consider $Q(\partial_t,\partial_x)$ an higher order partial differential operator defined as in \eqref{eq:Q_case1} or \eqref{eq:Q_case2}; in particular, suppose that $\mathcal{P}_m(\lambda,\xi)$ and $\mathcal{P}_{m-1}(\lambda,\xi)$ satisfy Hypothesis \ref{hyp:1} if  $Q=Q_1$, whereas  $\mathcal{P}_m(\lambda,\xi)$, $\mathcal{P}_{m-1}(\lambda,\xi)$ and $\mathcal{P}_{m-2}(\lambda,\xi)$ verify Hypothesis \ref{hyp:2} if $Q=Q_2$; consider the non-linear Cauchy problem 
\begin{equation}
\label{eq:CP_higher_order_nonlin}
\begin{cases}
Q(\partial_t,\partial_x)=|u|^{p_c}\mu(|u|),& \quad (t,x)\in [0,\infty)\times \R^n,\\
\partial_t^j u(0,x)=0,& \quad j\in \{0,\dots, m-1\}, \\
\partial_t^{m-1}u(0,x)=f(x),
\end{cases}
\end{equation}
in space dimension $n>m-2-\varrho$ where $\varrho=0$ if $Q=Q_1$, whereas $\varrho=1$ if $Q=Q_2$; here, $\mu:[0,\infty) \to [0,\infty)$ satisfies assumptions (i)-(iv), and $p_c$ is defined by
\begin{equation} 
\label{eq:pc_higher}
p_c= 1+\frac{m-\varrho}{n-(m-2-\varrho)}.
\end{equation}
\\
Assume that $f\in L^1(\R^n)$ satisfies the sign condition 
\begin{equation*}
\label{eq:sign_f_sigmaevol_ell=1}
 \int_{\R^n} f(x)\,dx>0;
\end{equation*}
if there exists a global weak solution $u\in L^{p_c}_{\text{loc}}\big([0,\infty),L^{p_c}_{\text{loc}}(\R^n)\big)$ to \eqref{eq:CP_higher_order_nonlin}, according to Definition \ref{def:weaksolution}, then
\begin{equation}
\label{eq:muint_higher}
\int_{0}^{c_0} \frac{\mu(s)}{s}\,ds<\infty,
\end{equation}
for some $c_0>0$, sufficiently small. 

On the other hand, suppose that the space dimension $n$ satisfies
\begin{equation} 
\label{eq:n_condtions}
m-2-\varrho<n< 2(m-1-\varrho);
\end{equation}
if $\mu=\mu(s)$ verifies the integral condition \eqref{eq:muint_higher}, then there exists $\eps_0>0$ such that for all $f\in L^1\cap L^2$ satisfying
\[\|f\|_{L^1\cap L^2}\leq \eps_0,\]
Cauchy problem \eqref{eq:CP_higher_order_nonlin} admits a unique global (in time) solution $u$ in the space $L^{\infty}_{loc}([0,\infty), L^{p_c}\cap L^\infty).$
\end{prop}

\begin{rem}
\label{rem:higher_order}
If $Q=Q_1$, the condition $\PI_{m-1}(0,\xi')=0$ in Hypothesis \ref{hyp:1} implies that $\lambda=0$ is a root of $\PI_{m-1}(\lambda,\xi')$ for all $\xi'\in S^{n-1}$, i.e. for all $\alpha\in \N^n$ with $|\alpha|=m-1$ it holds $c_{0,\alpha}=0$; in particular, since $\PI_{m-1}(\lambda,\xi)$ is strictly hyperbolic it is a simple root; as a consequence, there exists $\alpha \in \N^n$ with $|\alpha_1|=m-2$ such that $c_{1,\alpha_1}\neq 0$; on the other, since $\PI_m(\lambda,\xi')$ and $\PI_{m-1}(\lambda,\xi')$ strictly interlace, there exists $\alpha_2\in \N^n$ with $|\alpha_2|=m$ such that $c_{0,\alpha_2}\neq 0$; indeed, this latter condition guarantees that $\lambda=0$ is not a root of $\PI_m(\lambda,\xi')$. \\
Similarly, if $Q=Q_2$ then, as a consequence of Hypothesis \ref{hyp:2}, for all $\alpha\in \N^n$ with $|\alpha|=m-2$ it holds $c_{0,\alpha}=0$; moreover, there exist two multi-indices $\tilde \alpha_1, \tilde\alpha_2\in \N^n$ such that $|\tilde\alpha_1|=m-3$, $|\tilde\alpha_2|=m-1$ and the coefficients $c_{1,\tilde\alpha_1}$, $c_{0,\tilde\alpha_2}$ are both different from $0$.
\end{rem}
Under Hypothesis \ref{hyp:1} if $Q=Q_1$ and Hypothesis \ref{hyp:2} if $Q=Q_2$ we are able to evalutate explicitly the critical exponent to the semi-linear Cauchy problem \eqref{eq:CP_higher_order_nonlin} applying formula \eqref{eq:pcritical}. \\
Indeed, if $Q=Q_1$ then Remark \ref{rem:higher_order} allows to conclude that $j=0$ and $j=1$ belong both to the set $J$ defined in \eqref{eq:J}; in particular, $r_0=m$ and $r_1=m-2$. Moreover, since $c_{m,0}=1$ and $c_{m-1,0}>0$ we can deduce that $m, m-1\in J$ with $r_m=r_{m-1}=0$. Additionally, for all $j\in J$ with $2\leq j\leq m-2$ it holds $r_j\geq m-j-1$. \\
Similarly, if $Q=Q_2$ satisfies the assumptions given in Hypothesis \ref{hyp:2}, then we know that $j=0$ and $j=1$ belong to the set $J$; in particular, we have $r_0=m-1$ and $r_1=m-3$. Furthermore, since $c_{0,m}=1$ and the coefficients $c_{0,m-1}, \,  c_{0,m-2}$ are strictily positive we derive $m, m-1, m-2 \in J$ and $r_m=r_{m-1}=r_{m-2}=0$; in addition, for any $j\in J$ with $2\leq j\leq m-3$ it holds $r_j \geq m-j-2$.\\
Based on these preliminaries, we can evaluate
\[ g(\eta)=\begin{cases}
\eta(m-1-\varrho) \quad & \text{ if } \eta<1, \\
\eta+m-2-\varrho \quad & \text{ if } \eta\in [1,2], \\
m-\varrho \quad & \text{ if } \eta>2;
\end{cases}\]
as consequence, if $n>m-2-\varrho$ then $p_c\in (1,\infty)$ assumes value
\begin{equation*} 
p_c=h(2)= 1+\frac{m-\varrho}{n-(m-2-\varrho)}.
\end{equation*}
Then, the non-existence result in Proposition \ref{prop:higherorder} is an immediate consequence of Theorem \ref{thm:non-existenceCauchy}.

In \cite{DA2021} the author investigates long time decay estimates for the solution to the linear Cauchy problem associated to \eqref{eq:CP_higher_order_nonlin} with $Q$ an higher order differential operator in the form \eqref{eq:Q_case1} or \eqref{eq:Q_case2}. 
The following lemma collects some useful estimates given in Theorem 1 of \cite{DA2021}.
%
%
\begin{lem}
\label{prop:linear_higherorder}
Consider the linear Cauchy problem
\begin{equation}
\label{eq:CP_higher_order_linear}
\begin{cases}
Q(\partial_t,\partial_x)u=0, &\quad (t,x)\in [0,\infty)\times \R^n \\
\partial_t^j u(0,x)=0,& \quad j\in \{0,\dots, m-1\}, \\
\partial_t^{m-1}u(0,x)=f(x),
\end{cases}
\end{equation}
with $f$ in $L^1\cap L^2$. \\Assume that the polynomials $\PI_{m-1}(\lambda,\xi)$ and $\PI_{m}(\lambda,\xi)$ verify Hypothesis \ref{hyp:1}. Then, the solution $u^\lin=u^\lin(t,x)$ to \eqref{eq:CP_higher_order_linear} with $Q=Q_1$ satisfies the long time decay estimate
\begin{equation}
\label{eq:lin_estimate_higher1}
\| u^\lin(t,\cdot)\|_{\dot{H}^s}\leq C (1+t)^{-\frac{n}{4}-\frac{s-(m-2)}{2}}\|f\|_{L^1\cap L^2},
\end{equation}
for any $m-2-n/2<s\leq m-1$, with $C>0$ independent of $f$ and $t$.

Suppose that $\PI_{m-1}(\lambda,\xi)$, $\PI_{m}(\lambda,\xi)$ and $\PI_{m-2}(\lambda,\xi)$ verify Hypothesis \ref{hyp:2}. Then, the solution $u^\lin=u^\lin(t,x)$ to \eqref{eq:CP_higher_order_linear} with $Q=Q_2$ satisfies the long time decay estimate
\begin{equation}
\label{eq:lin_estimate_higher2}
\|u^\lin(t,\cdot)\|_{\dot{H}^s}\leq C (1+t)^{-\frac{n}{4}-\frac{s-(m-3)}{2}}\|f\|_{L^1\cap L^2},
\end{equation}
for any $m-3-n/2<s\leq m-1$, with $C>0$ independent of $f$ and $t$.
\end{lem}
Assumption \eqref{eq:n_condtions} on the space dimension $n$ guarantees that the critical exponent \eqref{eq:pc_higher} belongs to $[2,\infty)$; moreover, for any $q\in [p_c,\infty]$, it holds
\begin{equation}
\label{eq:n_conditions_consequences}
m-2-\frac{n}{2}-\varrho<n\Big(\frac{1}{2}-\frac{1}{q}\Big) <m-1;
\end{equation}
thus, as a consequence of Lemma \ref{prop:linear_higherorder}, using that $\|\cdot\|_{L^q}\lesssim \|\cdot\|_{\dot{H}^{n(\frac{1}{2}-\frac{1}{q})}}$ for any $q\in [2,\infty)$, we obtain the following decay estimate for the solution to the linear Cauchy problem \eqref{eq:CP_higher_order_linear},
\begin{equation*}
\|u^\lin(t,\cdot)\|_{L^q}\leq C (1+t)^{-\frac{n}{2}(1-\frac{1}{q})+\frac{m-2-\varrho}{2}}\|f\|_{L^1\cap L^2},
\end{equation*}
for any $q\in [p_c,\infty)$.
In particular, for any $q\in [p_c, \infty)$ we may estimate 
\begin{equation}
\label{eq:Lq_higher}
\| u^\lin(t,\cdot)\|_{L^q}\leq C (1+t)^{-\frac{1}{p_c}}\|f\|_{L^1\cap L^2};
\end{equation}
in order to estimate $\| u^\lin(t,\cdot)\|_ {L^\infty}$ for all $t\geq 0$ we employ the following lemma. 
\begin{lem}
Let $0<2s_1<n<2s_2$. Then, for any function $f\in \dot{H}^{s_1}\cap \dot{H}^{s_2}$ one has
\[ \|f\|_{\infty}\lesssim \|f\|_{\dot{H}^{s_1}}+\|f\|_{\dot{H}^{s_2}}.\]
\end{lem}
\noindent Due to \eqref{eq:n_conditions_consequences} we can consider $\eps>0$ sufficiently small such that
\[m-2-\frac{n}{2}-\rho+\eps<\frac{n}{2}\leq m-1-\eps,\]
and then estimates \eqref{eq:lin_estimate_higher1} and \eqref{eq:lin_estimate_higher2} holds true for $s=s_1:= n/2-\eps$ and $s=s_2:=n/2+\eps$. Thus, we may conclude that $ u^\lin(t,\cdot)\in L^\infty$ and, in particular, estimate \eqref{eq:Lq_higher} also applies to $q=\infty$. 

Since estimate \eqref{eq:Lq_higher} holds true for any $q\in [p_c,\infty]$, we can apply Theorem \ref{thm:existence} to conclude that under the assumptions of Proposition \ref{prop:higherorder} the non-linear Cauchy problem \eqref{eq:CP_higher_order_nonlin} admits a global (in time) small data solution provided that the integral condition \eqref{eq:muint_higher} holds true for some $c_0>0$.
\begin{ex}
As an example, one can consider the following higher order partial differential equation
\begin{equation*}
\big(\partial_t+\sigma \big)\big(\big(\partial_{t}^2-\mu\Delta\big)\big(\partial_t^2+\sigma\partial_t-c^2\Delta) -\gamma^2\partial_t^2\Delta) u=0,
\end{equation*}
which is related to a couple system of elastic waves with Maxwell equations in $\R^3$ (see Example 3.3. in \cite{DA2021}); here, $\mu>0$ is the Lam\`e constant, $\gamma>0$ is the coupling constant, $\sigma>0$ is the electric conductivity and $c^2 >0$ is the inverse of the product of the electric permittivity and magnetic permeability.
In this case the differential operator $Q(\partial_t,\partial_x)$ can be written as \[ Q(\partial_t,\partial_x)=\mathcal{P}_5(\partial_t,\partial_x)+\mathcal{P}_4(\partial_t,\partial_x)+\mathcal{P}_3(\partial_t,\partial_x),\]
where $\mathcal{P}_5(\partial_t,\partial_x)$, $\mathcal{P}_4(\partial_t,\partial_x)$ and $\mathcal{P}_3(\partial_t,\partial_x)$ have symbols respectively
\begin{align*}
\mathcal{P}_5(\lambda,\xi)&= \lambda(\lambda^4-(\mu+c^2+\gamma^2)\lambda^2|\xi|^2+c^2\mu|\xi|^4),\\
\mathcal{P}_4(\lambda,\xi)&= \sigma( 2\lambda^4-(2\mu+c^2+\gamma^2)\lambda^2|\xi|^2+c^2\mu|\xi|^4),\\
\mathcal{P}_3(\lambda,\xi)&= \sigma^2\lambda(\lambda^2-\mu|\xi|^2).
\end{align*}
It is possible to prove that $\mathcal{P}_5(\lambda, \xi')$, $\mathcal{P}_4(\lambda,\xi')$ and $\mathcal{P}_3(\lambda,\xi')$ satisfy Hypothesis \ref{hyp:2}; moreover, the space dimension $n=3$ satisfies the conditon \eqref{eq:n_condtions}; then, we can conclude that the critical exponent for the corresponding Cauchy problem with power non-linearity $|u|^p$ is $p_c=5$; in particular, if $f\in L^1\cap L^2$ is sufficiently small and verifies the sign condition 
\[\int_{\R^n}f(x)\,dx>0\]
then, the non-linear Cauchy problem
\begin{equation*}
\begin{cases}
&\big(\partial_t+\sigma \big)\big(\big(\partial_{t}^2-\mu\Delta\big)\big(\partial_t^2+\sigma\partial_t-c^2\Delta) -\gamma^2\partial_t^2\Delta) u=| u |^{5}\mu(|u|)\\
&\partial_t^j u(0,x)=0, \quad j\in \{0,\dots, 3\}, \\
&\partial_t^{4}u(0,x)=f(x),
\end{cases}
\end{equation*}
with $\mu$ satisfying properties (i)-(iv), admits a global (in time) solution in $L^{\infty}_{loc}([0,\infty), L^{5}\cap L^\infty)$ if, and only if, $\mu$ satisfies the integral condition \eqref{eq:muint_higher} for some $c_0>0$.

\end{ex}

\subsection{Damped $\sigma$-evolution equation with mass}
\label{sec:KG}
Applying Proposition \ref{thm:existence2} one can prove the following result about the existence of global (in time) solutions to the non-linear Cauchy problem for a damped $\sigma$-evolution equation with mass.
\begin{prop}
Consider for $n\geq 1$ the Cauchy problem
\begin{equation}
	\label{eq:dampedKG}
	\begin{cases}
		\partial_{t}^2 u+ (-\Delta)^\sigma u +  2a\partial_tu+m^2 u =|u|\mu(|u|),\\
		u(0,x)=0, \quad u_t(0,x)=f(x),
	\end{cases}
\end{equation}
for $(t,x)\in [0,+\infty)\times \R^n$, with $a>0$, $m>0$ and $\sigma>n/2$; assume that $\mu:[0,\infty)\to [0,\infty)$ is a continuous function such that $\mu(0)=0$ and  properties (i), (ii) and (iv) with $p_c=1$ hold true. Then, there exists $\eps_0>0$ sufficiently small such that for all $f\in L^2$ satisfying 
\[ \|f\|_{L^2}\leq \eps_0,\]
the Cauchy problem \eqref{eq:dampedKG} admits a unique global (in time) solution $u$ in $L^\infty([0,\infty), L^2\cap L^\infty)$.
\end{prop}
Employing the results obtained in \cite{DAmass}, being $\sigma>n/2$ we can derive for any $q\in [2,\infty]$ the following estimate for the solution $u^\lin=u^\lin(t,x)$ to the linear Cauchy problem associated to \eqref{eq:dampedKG},
\begin{equation*}
\|u^\lin(t,\cdot)\|_{L^q}\leq C e^{-c_1 t } \|f\|_{L^2},
\end{equation*}
where $C$ and $c_1$ are positive constants independent of $f$ and $t$.
In particular, for any $q\in [2,\infty]$ estimate \eqref{eq:Letap,Linfty_estimate_thm2} holds true for $\beta>1$ arbitrarily large, and $\delta=r=2$. As a consequence, if $f$ is sufficiently small in $ L^2$ the application of Proposition \ref{thm:existence2} allows to conclude that problem \eqref{eq:dampedKG} admits a unique global (in time) solution $u$ in $L^\infty([0,\infty), L^2\cap L^\infty)$ (see also Example \ref{example:F_pc=1}). \\
We mention that for $\sigma=1$ problem \eqref{eq:dampedKG} is known as damped Klein-Gordon equation.

\section{Concluding Remarks and Future Works}
\label{sec:conclusions}
Many authors have investigated long time decay estimates for the solution to linear Cauchy problems in the form
\begin{equation}
	\label{eq:CPlinear_conclusions}
	\begin{cases}
		L(\partial_t, \partial_x) u = 0, \quad \\
		\partial_t^j u(0,x)=0, \quad  \text{ for all } j=0,\dots, m-2,\\
		\partial_t^{m-1} u(0,x)= f(x).
	\end{cases}
\end{equation}
with $L$ as in \eqref{eq:Ldef}. For several models, given $\delta, r\in [1,\infty]$ properly chosen, it is possible to determine for all $q>\bar{q}\geq 1$ sufficiently large, an exponent $\beta^\delta_q$ depending on $\delta$ and $q$ such that for every $f\in L^\delta\cap L^r$ the solution to \eqref{eq:CPlinear_conclusions} satisfies
\begin{equation}
 \label{eq:decay_remark}
 \|u^\lin(t,\cdot)\|_{L^q}\leq C (1+t)^{-\beta^\delta_q}\|f\|_{L^\delta\cap L^r},  
\end{equation}
for some $C>0$ independent of $f$ and $t$; often, the exponent $\beta^\delta_q$ is strictly increasing with respect to $q\in [\bar{q},\infty]$. As already discussed, the obtained estimates may be applied to obtain some results about the existence of global (in time) small data solutions to the corresponding non-linear problem 
\begin{equation}
	\label{eq:CPnonlinear_conclusions}
	\begin{cases}
		L(\partial_t, \partial_x) = F(u), \quad \\
		\partial_t^j u(0,x)=0, \quad  \text{ for all } j=0,\dots, m-2,\\
		\partial_t^{m-1} u(0,x)= f(x),
	\end{cases}
\end{equation}
under suitable assumptions on the non-linear function $F$.
\begin{rem}
Assume $\delta=1$ and let $p_c>1$ the exponent obtained applying formula \eqref{eq:pcritical} with $\ell=0$; suppose that for $q=p_c$ the exponent $\beta^1_q$ in estimate \eqref{eq:decay_remark} satisfies $\beta_{p_c}^1=1/p_c$; then for any $f\in L^1\cap L^r$, sufficiently small, and $p>p_c$ a global (in time) solution to \eqref{eq:CPnonlinear_conclusions} with $F(u)=|u|^p$ exists; in this case, one can trivially conclude that a global (in time) solution still exists if one replace the power non-linearity $F(u)=|u|^p$ by the non-linear term $F(u)=|u|^p\mu(|u|)$ with $\mu:[0,\infty)\to [0,\infty)$ bounded and non-decreasing on $[0,\epsb]$, satisfying 
\[ F(0)=0,\qquad |F(y)-F(z)|\leq C\,|y-z|\,(|y|^{p-1}+|z|^{p-1})\,\mu(|y|+|z|),\]
for all $|y|, |z|\leq \epsb$, for some $\epsb>0$ arbitrarily small. \\ 
On the other hand, the critical exponent $p=p_c$ belongs to the non-existence range;  however, Theorem \ref{thm:existence} guarantees that Cauchy problem \eqref{eq:CPnonlinear_conclusions} with $F(u)=|u|^{p_c}\mu(|u|)$  still admits a global (in time) small data solution provided that $\mu$ satisfies additionally the  integral condition \eqref{eq:muint}.
Indeed, assumption \eqref{eq:muint} allows us to obtain
\begin{equation*}
\int_0^\frac{t}{2}(1+s)^{-1}\,\mu(c(1+s)^{-a})\,ds \leq  \int_1^\infty s^{-1}\,\mu(cs^{-a})\,ds = \frac1a\,\int_0^c \tau^{-1}\,\mu(\tau)\,d\tau,
\end{equation*}
for any $a>0$, by using the change of variable~$1+s\mapsto s$ first and $\tau=cs^{-a}$ later. As shown in the proof of Theorem \ref{thm:existence}, this latter estimate is crucial to prove the contraction argument which leads to the existence of the global small data solution.  
\end{rem}
\begin{rem}
If one drops the assumption of initial data to be (small) in $L^1$, it is well- known that, in general, the critical exponent for equation \eqref{eq:CPnonlinear_conclusions} with power non-linearity $|u|^p$ changes. In particular, suppose that the $L^1$ smallness is replaced by the $L^\delta$ smallness for some $\delta>1$. In this case a new critical exponent can be often determined if there exists  $\bar{p}_c>1$ such that  $\bar p_c>\bar{q}/\delta$ and the identity $\beta^\delta_{\delta \bar{p}_c} = 1/{\bar p_c}$ holds true; in particular, interestingly, for several models it is possible to prove that the non-linear problem \eqref{eq:CPnonlinear_conclusions} with $F(u)=|u|^{\bar p_c}$ admits a global (in time) small data solution (see, for instance, \cite{Ikeda2002, Ikehata2002, Ikehata2004 });
 this happens, in particular, if there exists $\gamma\in [1,\delta)$, $\gamma>\bar{q}/\bar{p}_c$ such that the solution to the linear problem \eqref{eq:CPlinear_conclusions} additionally satisfies
\begin{equation*}
 \label{eq:decay_remark3}
 \|u^\lin(t,\cdot)\|_{L^q}\leq C (1+t)^{-\beta^\gamma_q}\|f\|_{L^\gamma\cap L^r}, \quad \text{ for every } q\in [\gamma \bar p_c,\infty],
\end{equation*}
and moroever, the inequality
\begin{equation}
\label{eq:beta_condition_thm3}
\beta_q^\gamma+\beta_{\gamma \bar{p}_c}^\delta \bar{p}_c -1\geq \min\Big\{\beta_q^\delta, \frac{1}{\bar p_c}\Big\},
\end{equation} 
holds true for all $q\geq \gamma \bar{p}_c$.
Indeed, in this case the existence of a global small data solution to \eqref{eq:CPnonlinear_conclusions} with $F(u)=|u|^{\bar{p}_c}$ can be proved following the same approach used in the proof of Theorem \ref{thm:existence}, repleacing the definition of $\|\cdot\|_X$ by 
\[ \|u\|_{X}=\sup_{s\in [0,\infty)}\Big(\sup_{q\in [\gamma \bar p_c, \delta \bar p_c)}(1+s)^{\beta^\delta_{q}}\|u(s,\cdot)\|_{L^q}+\sup_{q\in [\delta \bar p_c,\infty]}(1+s)^{\frac{1}{\bar{p}_c}}\|u(s,\cdot)\|_{L^q}\Big),\]
and estimate \eqref{eq:Nu-Nv_estimate} by 
\begin{align*}
		\|Nu(t,\cdot)-Nv(t,\cdot)\|_{L^q}&\lesssim \int_0^\frac{t}{2} (1+t-s)^{-\beta_q^\gamma}\||u(s,\cdot)|^{\bar p_c}-|v(s,\cdot)|^{\bar{p}_c}\|_{L^\gamma\cap L^r}\,ds \\
		&+ \int_\frac{t}{2}^t (1+t-s)^{-\beta_q^\delta}\||u(s,\cdot)|^{\bar p_c}-|v(s,\cdot)|^{\bar{p}_c}\|_{L^\delta\cap L^r}\,ds.
	\end{align*} 
In particular, if the inequality \eqref{eq:beta_condition_thm3} holds true, then we may estimate
\begin{align*}
		\label{eq:int_0_t/2_1}
		\int_0^{\frac{t}{2}} \dots ds &\lesssim (1+t/2)^{-\beta_q^\gamma} \int_0^{\frac{t}{2}} (1+s)^{-\beta_{\gamma \bar p_c}^\delta  \bar{p}_c}\,ds \\
		&\lesssim (1+t)^{-\beta_q^\gamma+1-\beta_{\gamma \bar p_c}^\delta\bar p_c}\\
		&\lesssim (1+t)^{-\min\{\beta_q^\delta, \frac{1}{\bar{p}_c}\}},
	\end{align*}
being $\beta_{\gamma \bar p_c}^\delta<1/\bar p_c$. On the other hand, we have 
	\begin{align*}
		\int_{\frac{t}{2}}^t \dots ds & \lesssim (1+t/2)^{-1} \int_\frac{t}{2}^t (1+t-s)^{-\beta_q^\delta}\,ds \\
		& \lesssim (1+t)^{-1+(1-\beta_q^\delta)_+}\ell(t)\\
		& \lesssim (1+t)^{-\min\{\beta_q^\delta, \frac{1}{\bar{p}_c}\}}.
		\end{align*}
where $\ell(t)=\log(e+t)$ if $\beta^\delta_q=1$ and $\ell(t)=1$ otherwise.
Of course, in this case one can trivially conclude that a global (in time) solution still exists if one replace the power non-linearity $F(u)=|u|^{\bar{p}_c}$ by the non-linear term $F(u)=|u|^{\bar{p}_c}\mu(|u|)$ with $\mu:[0,\infty)\to [0,\infty)$ bounded and non-decreasing on $[0,\epsb]$, satisfying 
\begin{equation} 
\label{eq:F_examples}
F(0)=0,\qquad |F(y)-F(z)|\leq C\,|y-z|\,(|y|^{p_c-1}+|z|^{p_c-1})\,\mu(|y|+|z|),
\end{equation}
for all $|y|, |z|\leq \epsb$, for some $C>0$ and $\epsb>0$ arbitrarily small. \\
On the other hand, it remains an open problem to determine whether there is a necessary condition on $\mu$ for the existence of global (in time) small data solutions and how this condition is related to the choice of the space $L^\delta$ to which the initial data $f$ belongs. 
\begin{ex}
Consider the non-linear Cauchy problem for a classical damped wave equation 
\begin{equation}
\label{eq:dampedwave_delta-reg}
\begin{cases}
\partial_{t}^2 u- \partial_{xx}u + \partial_tu =|u|^p\mu(|u|), \quad (t,x)\in [0,+\infty)\times \R,\\
u(0,x)=0, \quad \partial_tu(0,x)=f(x),
\end{cases}
\end{equation}
with $f\in L^\delta \cap L^2$ for some $\delta \in (1,2]$. For any choice of $\delta$ in $(1,2]$ the solution to the linear Cauchy problem associated to \eqref{eq:dampedwave_delta-reg} satisfies the following decay estimate
\[ \|u^\lin(t,\cdot)\|_{L^q}\lesssim (1+t)^{-\beta_q^\delta}\|f\|_{L^\delta\cap L^2}, \quad \text{with } \beta_q^\delta=\frac{1}{2}\Big(\frac{1}{\delta}-\frac{1}{q}\Big),\]
for any $q\in [2,\infty]$. In particular, $\bar{p}_c=1+2\delta$ satisfies $\beta^\delta_{\delta \bar{p}_c}=1/\bar{p}_c$; moreover, for any $\gamma\in (1, \delta)$ it holds $\gamma>2/\bar{p}_c$ and 
\[ \beta_q^\gamma+\beta_{\gamma \bar{p}_c}^\delta \bar{p}_c -1=\beta_q^\delta,\]
for any $q\in [2,\infty]$; thus, condition \eqref{eq:beta_condition_thm3} is satisfied. \\ 
We can conclude that for all $p\geq \bar{p}_c=1+2\delta$ the Cauchy problem \eqref{eq:dampedwave_delta-reg} admits a global (in time) solution, provided that $\mu:[0,\infty)\to [0,\infty)$ is bounded and non-decreasing on $[0,\epsb]$ and satisfies inequality \eqref{eq:F_examples} for all $|y|,|z|\leq \bar{\eps}$, and $f$ is sufficiently small in $L^\delta\cap L^2$; on the other hand, if $\mu\equiv 1$ and $f$ satisfies $f(x)\geq \epsilon (1+|x|)^{-\frac{n}{\delta}}(\log(1+|x|))^{-1}$ for some $\epsilon>0$, then no global (in time) weak solutions to \eqref{eq:dampedwave_delta-reg} exist for all $1<p<1+2\delta$  (see, for instance, \cite{DAE2017}).\\
In the case $p=1+2\delta$ it remains an open problem to determine whether there is a necessary condition on $\mu$ for the existence of global (in time) small data solutions to \eqref{eq:dampedwave_delta-reg}.
\end{ex}
\end{rem}
\begin{rem}
If  $\delta>\bar{q}$ and $\beta^\delta_{\delta}>1$ in estimate \eqref{eq:decay_remark}, then assuming the initial data $f$ to be sufficiently small in $ L^\delta\cap L^r$, the non-linear problem \eqref{eq:CPnonlinear_conclusions} with $F(u)=|u|^p$ admits a global (in time) solution for any $p>1$. In general, it is not known if any global small data solution exists when $F(u)=|u|$; however, as Proposition \ref{thm:existence2} emphasizes, one can still prove the existence of global small data solutions to \eqref{eq:CPnonlinear_conclusions} if one replace, for instance,  $F(u)=|u|$ by the non-linear term $F(u)=|u|\mu(|u|)$ with $\mu:[0,\infty)\to [0,\infty)$ a continuous function satisfying $\mu(0)=0$ together with properties (i), (ii) and (iv) with $p_c=1$. On the other hand, it remains a challenging open problem to determine a necessary condition on $\mu$ for the existence of global (in time) small data solutions to problem \eqref{eq:CPnon-linear} with $F(u)=|u|\mu(|u|)$.
\end{rem}
\begin{rem}
In the present paper, we have considered evolution equations with constant coefficients. The case of differential operators $L$ with time-dependent coefficients will be the subject of future studies. It is well known that the presence of time-dependent coefficients can deeply influence the asymtpotic behavior of the solution to the linear problem, and then the critical exponent; in particular, formula \eqref{eq:pcritical} can not be applied in general to determine the critical power. Let us consider, for instance, the following Cauchy problem for a damped Klein-Gordon wave equation with time-dependent damping and mass terms
\begin{equation}
\label{eq:dampedKG_time_dep}
\begin{cases}
\partial_t^2 u -\Delta u + \lambda (1+t)^k \partial_t u +\nu^2 (1+t)^{2\ell}u=F(u),\\
u(0,x)=0, \quad \partial_t u(0,x)=f(x) 
\end{cases}
\end{equation}
for $(t,x)\in [0,\infty)\times \R^n$, with $\lambda, \nu>0$, $k\in (-1,1)$ and $\ell<k$. In \cite{DAGR2019} we considered $F(u)=|u|^p$ and we proved that the critical exponent is still the classical Fujita exponent $p_c=1+2/n$ if $\ell < (k-1)/2$; whereas, if $\ell = (k-1)/2$ then the critical power is given by 
\[ p_c= 1+\frac{2}{n+2\beta},\]
where $\beta=\nu^2/(\lambda(1-k))$. If finally $\ell>(k-1)/2$, then the Cauchy problem associated to \eqref{eq:dampedKG_time_dep} admits a global (in time) solution for any $p>1$, assuming $f$ to be small in $L^1\cap L^2$; the same holds true if the mass term dominates the damping term, namely $\ell> k$ or $\ell=k$ and $\nu>\lambda/4$ (see \cite{Girardi2019}).
In a future work we are going to investigate the existence of global (in time) solutions to \eqref{eq:dampedKG_time_dep} with non-linear term $F(u)=|u|^{p_c}\mu(|u|)$ in place of the power non-linearity $|u|^p$; our interest is to determine whether and how the presence of time-dependent coefficients influences the integral condition \eqref{eq:muint} that allows to precisely distinguish the region of existence of a global (in time) small data solution from that in which no global weak solution exists. Clearly, the presence of time-dependent coefficients introduces new difficulties, for instance, due to the fact that the solution to the linear problem is no longer invariant by translation. Moreover, it is necessary to employ a modified test function method (see \cite{DAL2013}) in order to investigate the non-existence of global (in time) weak solutions. \\
In \cite{DR} the authors consider model \eqref{eq:dampedKG_time_dep} in space dimension $n=1$, with $\nu=0$, $k\in [0,1)$ and $F(u)=|u|^3 \mu(|u|)$ with $\mu$ as in Example \ref{example:mu}; they prove that if $k\in [0,1)$ and $\lambda>0$ then the integral condition \eqref{eq:muint} still allows to describe a threshold between the existence and non-existence of global (in time) solutions. 
\end{rem}
\begin{rem}
For some choice of $L$ the critical exponent for problem \eqref{eq:CPnon-linear_classic} with power non-linearity, is not related to the scaling properties of the differential operator $L$ and then it cannot be evaluated applying formula \eqref{eq:pcritical}. This is the case, for instance, of the classical wave equation. It is well known that the semi-linear Cauchy problem
\begin{equation}
\label{eq:wave}
\begin{cases}
\partial_{tt} u- \Delta u =|u|^p, &\quad (t,x)\in [0,\infty)\times \R^n \\
u(0,x)=0, \quad \partial_t u(0,x)=f(x)&
\end{cases}
\end{equation}
with $f\in C_c^\infty(\R^n)$ admits a global (in time) small data solution if, and only if, $p>p_0(n)$ where $p_0(n)$ is the Strauss exponent defined as the positive root of the quadratic equation 
\[2+(n+1)p-(n-1)p^2=0.\]
%
For this model the integral condition \eqref{eq:muint} seems to be not suitable to identify a critical non-linearity; indeed, in the recent paper \cite{ChenReissig2024} the authors considered Cauchy problem \eqref{eq:wave} in space dimension $n=3$ with non-linearity $F(u)=|u|^{p_0(3)}\mu(|u|)$ in place of $|u|^p$; they provide a sharp condition on the behavior of $\mu(\tau)$ for $\tau\to 0^+$ which determine a threshold between global (in time) existence of small data radial solutions and blow-up of solutions in finite time; such condition is strictly related to the exponent $p_0(3)$.
\end{rem}
\begin{rem}
Given $L$ a partial differential operator, one can consider the Cauchy problem associated to \eqref{eq:CPnon-linear_classic} with a memory-type non-linearity in the form 
\[ F(\partial_t^\ell u)=\int_0^t (t-\tau)^{-\gamma}|\partial_t^\ell u(\tau,x)|^p d\tau,\]
with $\gamma\in (0,1)$ and $p>1$.
For this kind of non-linearities the critical exponent cannot in general be obtained by scaling arguments 
(see for instance \cite{Caz2008, DA2014NLA, DAG2021, Fino2011,Girardi2023}). In a forthcoming paper, we will replace $F$ by
\[ F(\partial_t^\ell u)=\int_0^t (t-\tau)^{-\gamma}|\partial_t^\ell u(\tau,x)|^p\mu(|\partial_t^\ell u(\tau,\cdot)|) d\tau,\]
for some $\mu:[0,\infty)\to [0,\infty)$, and we will investigates the influence of $\mu=\mu(s)$ on the global (in time) existence of small data solutions to the corresponding non-linear problem.
\end{rem}

\subsection*{Acknowledgment}
The author has been partially supported by his own INdAM-GNAMPA Project, Grant Code CUP E53C23001670001.

\end{document}